\documentclass{article}
\usepackage{amsmath,amsfonts,color}

\def\C{\mathbb{C}}

\def\pas{\par\smallskip}

\def\ba{\begin{array}}
\def\ea{\end{array}}

\def\zx#1{\begin{equation}\label{#1}}
\def\zc{\end{equation}}

\renewcommand{\theequation}{\thesection.\theequation}
\numberwithin{equation}{section}

\title{Degenerate 4-Dimensional Matrices with Semi-Group Structure}

    \author{ Olga Veko\footnote{Kalinkovichi Gymnasium, Belarus,vekoolga@mail.ru},
    Elena Ovsiyuk\footnote{Mosyr State Pedagogical University,
    Belarus, e.ovsiyuk@mail.ru},
 Alexandru Oana\footnote{Transilvania University of Brasov, alexandru.oana@unitbv.ro},
             Mircea Neagu\footnote{Transilvania University of
            Brasov, mircea.neagu@unitbv.ro},\\
            Vladimir Balan\footnote{University Politehnica of Bucharest, Romania, vladimir.balan@upb.ro},
             Victor Red'kov\footnote{B.I. Stepanov Institute
of Physics, NAS of Belarus, redkov@dragon.bas-net.by}}
    \date{}

 \date{}
\begin{document}\maketitle
\begin{abstract}

While dealing with the nontrivial task of classifying Mueller
matrices, of special interest is
    \index{degenerate!Mueller matrices} the study of the degenerate Mueller
    matrices (matrices with vanishing determinant, for which the law of multiplication holds, but
    there exists no inverse elements). Earlier, it was developed a special technique of parameterizing
    arbitrary 4-dimensional matrices with the use of a four 4-dimensional vector $(k,m,l,n)$. In the
    following, a classification of degenerate 4-dimensional real matrices of rank 1, 2, and 3 is elaborated.
To separate possible classes of degenerate matrices of ranks 1 and
2, we impose linear restrictions on
    $(k,m, l,n)$, which are compatible with the group multiplication law.
All the subsets of matrices obtained by this method, form either
subgroups or \index{semigroups} semi-groups. To obtain singular
matrices of rank 3, we specify 16 independent possibilities to get
the 4-dimensional
    matrices with zero
    determinant

\end{abstract}\pas\noindent


\section{introduction}

In polarization optics, an important issue is to classify possible
classes of the Mueller matrices.
    In particular, of special interest are degenerate Mueller matrices (with vanishing determinant).

There is known a special technique of parameterizing arbitrary
4-dimensional matrices from
    \index{gl@$GL(4,\C)$} $GL(4,\C)$ with the use of four 4-dimensional vector $ (k, m, l, n)$.
In the following, a classification of degenerate 4-dimensional
real matrices of rank 1, 2, 3 is elaborated. To separate possible
classes of \index{degenerate!matrices} degenerate matrices of
ranks 1 and 2, we
    impose linear restrictions on $ (k, m, l, n)$, which are compatible with the group multiplication law.
All the subsets of matrices obtained by this method  form either
subgroups or semigroups.
 To obtain singular matrices of rank 3, we specify 16 independent possibilities to get the 4-dimensional
 matrices with zero determinant\footnote{This approach  is based on the works
\cite{2008-Red'kov-Bogush-Tokarevskaya,
 Ovsiyuk-Veko-Red'kov-2012,
Ovsiyuk-Red'kov-Fert-2012,Ovsiyuk-Red'kov-RUDN,
Red'kov-Bogush-Tokarevskaya-2007}.}.

 In \index{Weyl spinor basis} Weyl spinor
basis, an arbitrary $(4\times 4)$-matrix can be parameterized
    by four 4-dimensional vectors $(k,m,l,n)$
\begin{eqnarray}G =\left(\ba{cc}
k_0 +\; {\bf k}\;\vec{\sigma}&n_0+\; {\bf n}\;\vec{\sigma}\\[3mm]
 l_0+\; {\bf l}\;\vec{\sigma}&m_0+\; {\bf m}\;\vec{\sigma}
\ea\right)=\left(\ba{cc}K&N\\[3mm]
 L&M \ea\right).\label{1}\end{eqnarray}
The matrices $G$ will be real if the second components of
parameters are imaginary
\begin{eqnarray}k_2^{*} =-k_2\;,\qquad m_2^{*} =-m_2\;,\qquad
    l_2^{*} =-l_2\;,\qquad n_2^{*} =-n_2\;,\nonumber\end{eqnarray}
and all remaining components are real. The law of multiplication
in explicit form is
\begin{eqnarray}k_0''=k_0'\; k_0+{\bf k}'\cdot {\bf k} +n'_0\; l_0 +{\bf n}' \cdot  {\bf l}\;,\nonumber\\
m_0''=m_0'\; m_0+{\bf m}'\cdot {\bf m}+l'_0\; n_0+{\bf l}'\cdot {\bf n}\;,\nonumber\\
n_0''=k_0'\; n_0+{\bf k}'\cdot {\bf n}+n'_0\; m_0+{\bf n}'\cdot
{\bf m}\;,\nonumber\\l_0''=l_0'\; k_0+{\bf l}'\cdot {\bf k}
+m'_0\; l_0+{\bf m}'\cdot {\bf l}\;,\nonumber\end{eqnarray}
\begin{eqnarray}{\bf k}''=k'_0\; {\bf k}+{\bf k}'\; k_0+i\; {\bf
k}'\times {\bf k}+n_0'\; {\bf l}+{\bf n}'\; l_0+i\;{\bf n}'\times {\bf l}\;,\qquad\nonumber\\
{\bf m}''=m'_0\; {\bf m}+{\bf m}'\; m_0+i\; {\bf
m}'\times {\bf m}+l_0'\; {\bf n}+{\bf l}'\; n_0+i\;{\bf l}'\times {\bf n}\;,\nonumber\\
{\bf n}''=k'_0\; {\bf n}+{\bf k}'\; n_0+i\; {\bf
k}'\times {\bf n}+n_0'\; {\bf m}+{\bf n}'\; m_0+i\;{\bf n}'\times {\bf m}\;,\nonumber\\
{\bf l}''=l'_0\; {\bf k}+{\bf l}'\; k_0+i\; {\bf l}'\times {\bf
k}+m_0'\; {\bf l}+{\bf m}'\; l_0+i\; {\bf m}'\times {\bf l
}\;.\label{2}\end{eqnarray}

Now we will study the degenerate 4-dimensional matrices of rank 1,
2, 3.  To separate possible classes
    of degenerate matrices of ranks 1 and 2, we impose linear restrictions on $ (k, m, l, n)$, which
    are compatible with the multiplication law (\ref{2}). All the subsets of matrices obtained by this method,
    form either sub-groups or semigroups. Below, only final results will be presented.
\section{One Independent Vector:  Variant I(k)}
First, we consider the variants with one independent vector.

\vspace{3mm}

{\bf Variant $I(k)$}:
\begin{eqnarray}{\bf n}=A\; {\bf k}\;,\; n_0=\alpha\; k_0\;,\;
{\bf m}=B\; {\bf k}\;,\; m_0=\beta\; k_0\;,\; {\bf l}=D\; {\bf
k}\;,\; l_0=t\; k_0\;.\label{B3.1x}\end{eqnarray}
The group multiplication law takes the form
\begin{eqnarray}\left.\ba{l} k_0''=k_0'\; k_0+{\bf k}'\cdot
{\bf k} +\alpha t\;
k'_0 k_0+A D\; {\bf k}'\cdot {\bf k}\;,\\[2mm]
m''_0=\beta^2\; k_0' k_0+B^2\; {\bf k}'\cdot {\bf k}
+t\alpha\; k'_0\; k_0+D A\; {\bf k}'\cdot {\bf k}\;,\\[2mm]
 n_0''=\alpha\; k_0' k_0+A\; {\bf k}'\cdot {\bf k}
+\alpha\beta\; k'_0\; k_0+A B\; {\bf k}'\cdot {\bf k}\;,\\[2mm]
l_0''=t\; k_0'\; k_0+D\; {\bf k}'\cdot {\bf k} +\beta t\; k'_0
k_0+BD\; {\bf k}'\cdot {\bf k}\;,
 \ea \right.
\nonumber\end{eqnarray}
\begin{eqnarray}\left.\ba{l} {\bf k}''=k'_0\; {\bf k}+{\bf k}'
\; k_0+i\; {\bf k}'\times {\bf k}+\alpha D\; k_0'\; {\bf k}+A t\;
{\bf k}' k_0+i A D\; {\bf k}'\times {\bf k}\;,\\[2mm]
 {\bf m}''=\beta B\; k'_0 {\bf k}+B\beta\; {\bf k}'\; k_0
+iB^2\; {\bf k}'\times {\bf k} + tA\; k_0' {\bf k}+D \alpha\; {\bf
k}'\; k_0 +
 i D A\; {\bf k}'\times {\bf k}\;,\\[2mm]
{\bf n}''=A\; k'_0\; {\bf k}+\alpha\; {\bf k}'\; k_0 +
 i A\; {\bf k}'\times {\bf k} +
\alpha B\; k_0' {\bf k}+A\beta\; {\bf k}'\; k_0+i AB
\; {\bf k}'\times {\bf k}\;,\\[2mm]
{\bf l}''=t\; k'_0 {\bf k}+D\; {\bf k}' k_0+iD\; {\bf k}'\times
{\bf k}+\beta D\; k_0'\; {\bf k}+B t\; {\bf k}'\; k_0+i B D\; {\bf
k}'\times {\bf k }\;. \ea \right. \nonumber\end{eqnarray}

Let us impose the relations (\ref{B3.1x}) for the double primed
parameters:

\vspace{2mm} $ \underline{n_0''=\alpha\; k_0''\;,}$
\begin{eqnarray}\alpha (1+\beta)\; k_0' k_0+A (1+B)\; {\bf k}'\cdot {\bf
k}=\alpha (1+\alpha t)\; k_0'\; k_0+\alpha (1+AD) \; {\bf k}'\cdot
{\bf k}\;,\nonumber\end{eqnarray}

$ \underline{ m_0''=\beta\; k_0''\;,} $
\begin{eqnarray}(\beta^2+t\alpha)\; k_0' k_0+(B^2 +DA)\; {\bf k}'\cdot {\bf k}=\beta (1+\alpha t)\;
    k_0'\; k_0+\beta(1+AD)\; {\bf k}'\cdot {\bf k}\;,\nonumber \end{eqnarray}

$ \underline{ l_0 ''= t\; k_0 ''\;, } $
\begin{eqnarray}t(1+\beta)\; k_0'\; k_0+D(1+B)\; {\bf k}'\cdot {\bf k}
= t\;(1+\alpha t)\; k_0'\; k_0+t (1+AD)\; {\bf k}'\cdot {\bf
k}\;,\nonumber\end{eqnarray}

$\underline{ {\bf n}''=A\; {\bf k}''\;, } $
\begin{eqnarray}(A+\alpha B)\; k'_0\; {\bf k}+(\alpha+A\beta)\; {\bf k}'\cdot k_0 +
 i A (1+B)\; {\bf k}'\times {\bf k}\nonumber\\
= A (1 +\alpha D)\; k'_0\; {\bf k}+A (1+A t)\; {\bf k}'\; k_0 +  i
A(1 +AD)\; {\bf k}'\times {\bf k}\;, \nonumber\end{eqnarray}

\vspace{5mm} $ \underline{{\bf m}''=B\; {\bf k}''\;, } $
\begin{eqnarray}(\beta B+tA)\; k'_0 {\bf k}+(B\beta+D \alpha )\;
{\bf k}'\; k_0 +i (B^2+AD)\; {\bf k}'\times {\bf k}
\nonumber\\
=B (1+\alpha D) k'_0\; {\bf k} + B (1+A t)\; {\bf k}'\; k_0+iB
(1+AD)\; {\bf k}' \times {\bf k}\;, \nonumber
 \end{eqnarray}

$\underline{ {\bf l}''=D\; {\bf k} ''\;,}$
\begin{eqnarray}(t+\beta D)\; k'_0 {\bf k}+(D +Bt)\; {\bf k}' k_0
+ iD (1+B)\; {\bf k}'\times {\bf k}\nonumber\\
=
 D (1+\alpha D)\; k'_0\; {\bf k}+D (1+ At)\;
{\bf k}'\; k_0 +
 i D(1+AD)\; {\bf k}'\times {\bf k}\;.
\nonumber\end{eqnarray}
We get the equations:
\begin{eqnarray} \alpha (1+\beta)=\alpha (1+\alpha t)\;,\qquad
 A (1+B)=\alpha (1+AD)\;,\nonumber\\
(\beta^2+t\alpha)=\beta (1+\alpha t)\;,\qquad (
B^2 +DA)=\beta(1+AD)\;,\nonumber\\
t (1+\beta)=t (1+\alpha t)\;,\qquad D (1 +B)=t(1 + AD)\;,
\nonumber\end{eqnarray}
\begin{eqnarray}(A+\alpha B)=A (1 +\alpha D)\;,\;
 (\alpha+A\beta)=A (1+A t)\;,\; A (1+B)=A(1 +AD)\;,\nonumber\\
 (\beta B+tA)= B (1+\alpha D)\;,\;(B\beta+D \alpha)=B (1+At)\;,\;
 (B^2+AD)= B (1+AD)\;,\nonumber\\
 (t+\beta D)=D (1+\alpha D)\;,\;(D +Bt)=D (1+ At)\;,\; D (1+B)=D(1+AD)\;.
 \nonumber\\\label{B3.2x}\end{eqnarray}

We are able now to find all the solutions for this system. First,
let us consider the
    system (\ref{B3.2x}) with three vanishing blocks -- we obtain only

\vspace{3mm}

\underline{solution $(K-1)$},
\begin{eqnarray}A=\alpha =0\;,\; B=\beta=0\;,\; D= t=0\;,\qquad
 G =\left(\ba{cc} K&0\\0&0 \ea\right) ;\label{B3.3x}\end{eqnarray}
these 4-matrices form a  semi-group. Their rank is equal to 2.
    When $\det\; K =0$, then the rank of $G$ equals to 1.

Let us consider the solutions for two vanishing blocks. There
arise three possibilities.

First, let it be
\begin{eqnarray}A=\alpha =0\;,\qquad D= t=0\; ;
\label{B3.4a}
 \end{eqnarray}
the system (\ref{B3.2x}) gives
\begin{eqnarray}\beta^2=\beta\;,\qquad B^2=\beta\;,\qquad\beta B=B\;,\qquad B^2=B\;,
\label{B3.4b}
 \end{eqnarray}
from which we obtain one new solution when $B=\beta =0$:

\vspace{3mm}

\underline{solution $(K-2)$},
\begin{eqnarray}A =\alpha =0\;,\; D= t=0\;,\; D=\beta=+1\;,\qquad
G=\left(\ba{cc}k_0+{\bf k}\vec{\sigma}&0\\
0&k_0+{\bf k}\vec{\sigma}\ea\right) ;\label{B3.4c} \end{eqnarray}
it is a set of non-degenerate matrices with a group structure.

Now, let it be
\begin{eqnarray}A =\alpha =0\;,\qquad B=\beta =0\; ;\label{B3.5a} \end{eqnarray}
the system (\ref{B3.2x}) gives
\begin{eqnarray}
0=0\;,\; A=0\;,\; 0=0\;,\; 0=0\;,\; t=t\;,\; D=t\;,\nonumber\\
0=0\;,\;  0=0\;,\;  0=0\;,\; 0=0\;,\;  0=0\;,\;  0=0\;,\; t=D\;,
\;  D =D\;,\;  D=D\;, \nonumber \label{B3.5b}\end{eqnarray}
so we obtain

\vspace{3mm} \underline{solution $(K-3)$},
\begin{eqnarray}A =\alpha =0\;,\; B=\beta =0\;,\; D=t\;,
G=\left(\ba{cc}K&0\\DK&0\ea\right),\nonumber\\
 G'G= \left(\ba{cc}K'&0\\DK'&0\ea\right)
\left(\ba{cc}K&0\\DK&0\ea\right)=\left(\ba{cc}
K'K&0\\DK'K&0\ea\right);\label{B3.5c}\end{eqnarray}
this is a set of degenerate matrices of the rank 2 with a
semi-group structure.

Now, let it be
\begin{eqnarray}B=\beta=0\;,\qquad D=t=0\; ;\label{B3.6a} \end{eqnarray}
the system (\ref{B3.2x}) reads
\begin{eqnarray} \alpha=\alpha\;,\qquad A=\alpha\;,\quad
 0=0\;,\qquad 0=0\;,\quad
 0=0\;,\qquad 0=0\;,\nonumber\\
 A=A\;,\;\; \alpha=A\;,\;\;  A=A\;,\;\; 0=0\;,\;\;  0=0\;, \;\;  0=0\;,\;\;
 0=0\;,\;\;  0=0\;,\;\; 0=0\;,
 \nonumber
 \label{B3.6b}\end{eqnarray}
so we get

\vspace{3mm}

\underline{solution $(K-4)$},
\begin{eqnarray}A= \alpha\;,\; B=\beta=0\;,\; D=t=0\;,\; G=\left(\ba{cc}K&A K\\
0&0\ea\right),\nonumber\\G'G=\left(\ba{cc}K'&A K'\\0&0\ea\right)
\left(\ba{cc}K&A K\\0&0\ea\right)=\left(\ba{cc}K'K&A K'K\\
0&0\ea\right);\label{B3.6c} \end{eqnarray}
it is a set of matrices of the rank 2 and with a semi-group
structure.

Let us consider cases with one vanishing block. The first
possibility is
\begin{eqnarray}A= \alpha=0\; ;\label{B3.7a}\end{eqnarray}
the system (\ref{B3.2x}) reads
\begin{eqnarray} 0=0,\; 0=0,\;\beta^2 =\beta,\; B^2=\beta,\;
t (1+\beta)=t,\; D (1 +B)=t,\;0=0,\; 0=0,\; 0=0\;,\nonumber\\
\beta B=B\;,\; B\beta=B\;,\; B^2=B\;,\; (t+\beta D)=D\;,\; (D
+Bt)=D\;,\;D (1+B)=D\;. \nonumber \label{B3.7b}\end{eqnarray}
These equations  lead to already known ones:

\vspace{3mm}

\underline{solution $(K-3)$},
\begin{eqnarray}A=\alpha=0\;,\; B=\beta=0\;,\; t= D\;,\; G=\left(
    \ba{cc} K&0\\D K&0\ea\right) ;\label{B3.7c}\end{eqnarray}

\underline{solution $(K-2)$},
\begin{eqnarray}A=\alpha=0\;,\; B=\beta=+1\;,\; D=t =0\;,\; G=\left(\ba{cc}k_0+{\bf k}
    \vec{\sigma}&0\\0&k_0+{\bf k}\vec{\sigma}\ea\right).\label{B3.7d} \end{eqnarray}

The second possibility is
\begin{eqnarray}D=t=0\; ;\label{B3.8a}\end{eqnarray}
the system (\ref{B3.2x}) leads to
\begin{eqnarray} \alpha\beta=0\;,\;A (1+B)=\alpha\;,\;
\beta^2 =\beta\;,\; B^2=\beta\;,\nonumber\\
0=0\;,\; 0=0\;,\;(A+\alpha B)=A\;,\; (\alpha+A\beta)=A\;,\;\; A
B=0\;,
\nonumber\\
\beta B=B\;,\; B\beta=B\;,\qquad B^2=B\;, \; 0=0\;,\; 0=0\;,\;
0=0\;. \nonumber \label{B3.8b}\end{eqnarray}
The equations (\ref{B3.8b}) leads to already known ones:

\underline{solution $(K-2)$},
\begin{eqnarray}A= \alpha =0\;,\; B=+1\;,\;\beta=+1\;,\; D=t =0\;,\nonumber\\
G=\left(\ba{cc}k_0+{\bf k}\vec{\sigma}&0\\
0&k_0+{\bf k}\vec{\sigma}\ea\right) ;\qquad\qquad\label{B3.8c}
\end{eqnarray}

\underline{solution $(K-1)$},
\begin{eqnarray}A=\alpha\;,\; B=\beta=0\;,\; D=t =0\;,\qquad G =\left(
    \ba{cc} K&0\\0&0 \ea\right).\label{B3.8d}\end{eqnarray}

The third possibility is
\begin{eqnarray}B=\beta=0\; ;\label{B3.9a}\end{eqnarray}
the system (\ref{B3.2x}) gives
\begin{eqnarray}
0=\alpha^2 t\;,\; A=\alpha (1+AD)\;,\; t\alpha=0\;,\; DA=0\;,\nonumber\\
0=\alpha t^2\;,\;D =t(1+AD)\;,\;0=A \alpha D\;,\;  \alpha=A (1+A t)\;,\; 0=A^2 D\;,\nonumber\\
tA=0\;,\; D \alpha=0\;,\; AD= 0\;,\; t=D (1+\alpha D)\;,\; 0 =D
At\;,\; 0= AD^2\;. \nonumber \label{B3.9b}\end{eqnarray}
Here there arise the already known solutions:

\vspace{2mm} \underline{solution $(K-4)$},
\begin{eqnarray}A=\alpha\;,\; B=\beta=0\;,\; D=t =0\;,\qquad G=\left(
    \ba{cc}K&A K\\0&0\ea\right) ;\label{B3.9c} \end{eqnarray}

\underline{solution $(K-3)$},
\begin{eqnarray}A =\alpha =0\;,\; B=\beta =0\;,\; D=t\;,\qquad G=\left(\ba{cc}
K&0\\DK&0\ea\right).\label{B3.9d} \end{eqnarray}

 Now, let us consider the cases in which all the blocks are non-vanishing
\begin{eqnarray}A\;,\; \alpha \neq 0\;,\qquad B\;,\;\beta \neq 0\;,\qquad
    D\;,\; t \neq 0\; ;\label{B3.10a}\end{eqnarray}
then the system (\ref{B3.2x}) gives
\begin{eqnarray}\beta =+\alpha t\;,\qquad  A (1+B)=\alpha (1+AD)\;,\nonumber\\
(\beta^2+t\alpha)=\beta (1+\alpha t)\;,\qquad (B^2 +DA)=\beta(1+AD)\;,\nonumber\\
\beta =+\alpha t\;,\qquad D (1 +B)=t(1+AD)\;,\nonumber\\
 B=A D\;,\qquad\;\; (\alpha+A\beta)=A (1+A t)\;,\qquad B=AD\;,\nonumber\\
(\beta B+tA)= B (1+\alpha D)\;,\; (B\beta+D \alpha)=B (1+A t)\;,\; (B^2+AD)= B (1+AD)\;,\nonumber\\
(t+\beta D)=D (1+\alpha D)\;,\;\;\;\qquad B=AD\;,\qquad B=
AD\;.\nonumber\\\label{B3.10b}\end{eqnarray}
By eliminating $B$ and $\beta$ by $B= AD\;,\;\beta=\alpha t$, the
remaining independent equations are
\begin{eqnarray} (A-\alpha)(1+AD) =0\;,\;\;
(AD-\alpha t)(1+AD) =0\;,\;\; (D -t)(1+AD) =0\;,
\nonumber\\
 (A- \alpha)(1+A t) =0\;,\quad (D-t)(1+\alpha D)=0\;.\qquad\qquad
 \label{B3.10c}\end{eqnarray}

First, let us study the case $A=\alpha$; the system (\ref{B3.10c})
takes the form
\begin{eqnarray}A= \alpha\;,\qquad (D -t)(1+AD) =0\;.\label{B3.11}\end{eqnarray}
We  get  two new solutions:

\vspace{3mm}

\underline{solution $(K-5)$},
\begin{eqnarray}A=\alpha\;,\; B=\beta=AD\;,\; D= t\;,\qquad G=\left(\ba{cc}
K&A K\\D K&AD K\ea\right) ;\label{B3.12a}\end{eqnarray}
it is a set of degenerate matrices of  rank 4 with semi-group
structure. As easily verified, the multiplication law holds
indeed:
\begin{eqnarray}\left(\ba{cc}K'&A K'\\D K'&AD K'\ea
\right) \left(\ba{cc}K&A K\\D K&AD K\ea\right)\nonumber\\
= \left(\ba{cc}(K' K+ ADK'K)&A (K' K+ ADK'K)\\
D (K' K+ ADK'K)&AD (K' K+
ADK'K)\ea\right).\label{B3.12b}\end{eqnarray}
%

\vspace{3mm}

\underline{Solution $(K-5)'$},
\begin{eqnarray}A=\alpha\;,\; B=\beta=-1\;,\; D= t\;,\; G=\left(\ba{cc}
K&A K\\-A^{-1} K&-K\ea\right),\; G'G=0\;
;\label{B3.12c}\end{eqnarray}

\underline{solution $(K-6)$},
\begin{eqnarray}A=\alpha\;,\; B=-1,\;\beta =- A t,\; D= -{1 \over A}\;,\nonumber\\
G=\left(\ba{cc} k_0+ {\bf k}\vec{\sigma}&Ak_0+ A {\bf k}\vec{\sigma}\\
tk_0-A^{-1} {\bf k}\vec{\sigma}&At k_0- {\bf k}
\vec{\sigma}\ea\right) ; \label{B3.13a}\end{eqnarray}
they are sets of degenerate matrices of rank 2 with semi-group
structure.

Let us verify the multiplication law in the case $K-6$:
\begin{eqnarray}G'G=\left(\ba{cc}k_0'+ {\bf k}\vec{\sigma}&Ak_0'+ A {\bf k}'\vec{\sigma}\\
    tk_0'-A^{-1} {\bf k}'\vec{\sigma}&At k_0'- {\bf k}'\vec{\sigma}\ea\right) \left(\ba{cc}
    k_0+ {\bf k}\vec{\sigma}&Ak_0+ A {\bf k}\vec{\sigma}\\
    tk_0-A^{-1} {\bf k}\vec{\sigma}&At k_0- {\bf k}\vec{\sigma}\ea\right)=G''\; ;\nonumber\end{eqnarray}
the result is presented via blocks:
\begin{eqnarray}(11)=(1+At)\; k'_0 k_0+(1+At) {\bf k}' k_0 \vec{\sigma}\;,\; \quad\;  \quad \;
\nonumber\\
(12)=A\; \left[ (1+At)\; k'_0k_0+(1+At) {\bf k}' k_0\vec{\sigma} \right],\nonumber \quad \\
(21)=t\; (1+At)\; k'_0 k_0-A^{-1}\; (1+At) {\bf k}' k_0\vec{\sigma}\;,\, \nonumber\\
(22)=At\; (1+At)\; k'_0 k_0-(1+At) {\bf k}' k_0
\vec{\sigma}\;,\quad \;\; \nonumber\end{eqnarray}
that is
\begin{eqnarray}k_0''=(1+At)k_0\; k'_0\;,\qquad {\bf k}''= (1+At)k_0\; {\bf k}'\;.\label{3.13b}\end{eqnarray}

Let us turn back to the system (\ref{B3.10c}) and consider the
case
\begin{eqnarray}A \neq \alpha\;,\qquad 1+AD=0\; ;\label{B3.14a}\end{eqnarray}
the equations (\ref{B3.10c}) give only

\vspace{3mm}

\underline{solution $(K-7)$},
\begin{eqnarray}A, \alpha,\; B= -1,\;\beta=-{\alpha \over A}\;,\; D= t=-{1 \over A}\,,\nonumber\\
 G=\left(\ba{cc}k_0+ {\bf k}\vec{\sigma}&\alpha k_0+ A {\bf k}\vec{\sigma}\\
-A^{-1} (k_0+{\bf k}\vec{\sigma})&-A^{-1} \alpha k_0 -{\bf
k}\vec{\sigma}\ea \right) ;\label{B3.14c}\end{eqnarray}
this is a set of degenerate matrices of rank 2 with semi-group
structure. Let us verify the multiplication law:
\begin{eqnarray}
\hspace{-3mm} G' G=\left(\ba{cc}
k_0'+ {\bf k}'\vec{\sigma}&\alpha k_0'+ A {\bf k}'\vec{\sigma}\\
-A^{-1} (k_0'+{\bf k}'\vec{\sigma})&-A^{-1} \alpha k_0'-{\bf
k}'\vec{\sigma} \ea\right) \left(\ba{cc}
k_0+ {\bf k}\vec{\sigma}&\alpha k_0+ A {\bf k}\vec{\sigma}\\
-A^{-1} (k_0+{\bf k}\vec{\sigma})&-A^{-1} \alpha k_0 -{\bf
k}\vec{\sigma} \ea\right);\nonumber\end{eqnarray}
the result is presented by blocks
\begin{eqnarray}
(11)=(1-{\alpha \over A})\; k_0' k_0+(1-{\alpha \over
A}) k'_0 {\bf k}\;\vec{\sigma}\;, \qquad \qquad \; \;\; \,\nonumber\\
(12)=\alpha\; (1-{\alpha \over A})\; k_0' k_0+A\; (1
- {\alpha \over A}) k'_0 {\bf k}\;\vec{\sigma}\;, \quad \quad \, \;\;\; \nonumber\\
(21)=-A^{-1}\; (1-{\alpha \over A})\; k_0' k_0 -A^{-1}
\; (1-{\alpha \over A}) k'_0 {\bf k}\;\vec{\sigma}\;,\nonumber\\
(22)=-A^{-1} \alpha\; (1-{\alpha \over A})\; k_0' k_0 - (1-{\alpha
\over A}) k'_0 {\bf k}\;\vec{\sigma}\;, \quad \;\;
\nonumber\end{eqnarray}
so we get
\begin{eqnarray}G''=G'G\;,\qquad k''_0=(1-{\alpha \over A})\; k_0'
k_0\;,\qquad {\bf k}''= (1-{\alpha \over A}) k'_0 {\bf k}
\;.\label{B3.14d}\end{eqnarray}

Thus, the analysis of the variant {\bf I(k)} is completed.
\section{One Independent Vector: Variant I(m)}
Let us consider {\bf variant $I(m)$}:
\begin{eqnarray}{\bf n}=A\; {\bf m}\;,\; n_0=\alpha\; m_0\;,\;
{\bf k}=B\; {\bf m}\;,\; k_0=\beta\; m_0\;,\; {\bf l}=D\; {\bf
m}\;,\;\;\; l_0=t\; m_0\;.\label{B4.1x}\end{eqnarray}
The group multiplication law takes the form
\begin{eqnarray}k_0''=\beta m_0'\;\beta m_0+B {\bf m}'\; B {\bf m}
+\alpha m'_0\; t m_0+A {\bf m}'\; D {\bf m}\;,\nonumber\\
m_0''=m_0'\; m_0+{\bf m}' \cdot  {\bf m}
+t m'_0\; \alpha m_0+D {\bf m}'\; A {\bf m}\;,\nonumber\\
n_0''=\beta m_0'\; \alpha m_0+B {\bf m}'\; A {\bf m}
+\alpha m'_0\; m_0+A {\bf m}'\cdot {\bf m}\;,\nonumber\\
l_0''=t m_0'\;\beta m_0+D {\bf m}'\; B {\bf m} +m'_0\; t m_0+{\bf
m}'\; D {\bf m}\;, \nonumber\label{B4.2a}\end{eqnarray}
\begin{eqnarray} {\bf k}''=\beta m'_0\; B {\bf m}+B {\bf m}'\;\beta m_0+i\; B {\bf
m}'\times B {\bf m}+\alpha m_0'\; D {\bf m}+A {\bf m}'\;
t m_0+i\; A {\bf m}'\times D {\bf m}\;,\nonumber\\
{\bf m}''=m'_0\; {\bf m}+{\bf m}'\; m_0+i\; {\bf m}'\times {\bf
m}+t m_0'\; A {\bf m}+D {\bf m}'\; \alpha
m_0+i\; D {\bf m}'\times A {\bf m}\;,\nonumber\\
{\bf n}''=\beta m'_0\;A {\bf m}+B{\bf m}'\; \alpha m_0 + i\; B
{\bf m}'\times A {\bf m}+\alpha m_0'\; {\bf m} +
A {\bf m}'\; m_0+i\; A {\bf m}'\times {\bf m}\;,\nonumber\\
 {\bf l}''=t m'_0\;B {\bf m}+D {\bf m}'\;\beta m_0+i\; D {\bf
m}'\times B {\bf m}+m_0'\;D {\bf m}+{\bf m}'\; t m_0 + i\; {\bf
m}'\times D {\bf m }\;. \nonumber
 \end{eqnarray}
We require
\begin{eqnarray}n_0''=\alpha\; m_0''\qquad \Longrightarrow\nonumber\\
\beta m_0'\; \alpha m_0+B {\bf m}'\; A {\bf m}
+\alpha m'_0\; m_0+A {\bf m}'\cdot {\bf m}\nonumber\\
=\alpha m_0'\; m_0+\alpha {\bf m}'\cdot {\bf m} +\alpha t m'_0\;
\alpha m_0+\alpha D {\bf m}'\; A {\bf m}\qquad
\Longrightarrow\nonumber\\\alpha (\beta+1)=\alpha (1+\alpha
t)\;,\qquad A(B+1)=\alpha(1 +AD)\; ; \nonumber\end{eqnarray}
\begin{eqnarray}k_0''=\beta\; m_0''\qquad \Longrightarrow\nonumber\\
\beta m_0'\;\beta m_0+B {\bf m}'\; B {\bf m} +\alpha
m'_0\; t m_0+A {\bf m}'\; D {\bf m}\nonumber\\
=\beta m_0'\; m_0+\beta {\bf m}'\cdot {\bf m}
+\beta t m'_0\; \alpha m_0+\beta D {\bf m}'\; A {\bf m}\qquad \Longrightarrow\nonumber\\
\beta^2+\alpha t=\beta (1+\alpha t)\;,\qquad B^2 +AD =\beta
(1+AD)\; ; \nonumber\end{eqnarray}
\begin{eqnarray}l_0''=t\; m_0''\qquad \Longrightarrow\nonumber\\
t m_0'\;\beta m_0+D {\bf m}'\; B {\bf m}
+m'_0\; t m_0+{\bf m}'\; D {\bf m}\nonumber\\
= t m_0'\; m_0+t {\bf m}'\cdot {\bf m}
+t^2 m'_0\; \alpha m_0+t D {\bf m}'\; A {\bf m}\qquad \Longrightarrow\nonumber\\
t (\beta +1)=t(1+\alpha t)\;,\qquad D(B +1)=t(1+AD) \; ;
\nonumber\end{eqnarray}
\begin{eqnarray}{\bf k}''=B\; {\bf m} ''\qquad \Longrightarrow\qquad\qquad\qquad\nonumber\\
\beta m'_0\; B {\bf m}+B {\bf m}'\;\beta m_0+i\; B {\bf m}'\times
B {\bf m}+\alpha m_0'\; D {\bf m}+A {\bf
m}'\; t m_0+i\; A {\bf m}'\times D {\bf m}\nonumber\\
= B m'_0\; {\bf m}+B {\bf m}'\; m_0+i\; B {\bf m}' \times {\bf
m}+B t m_0'\; A {\bf m}+B D {\bf m}'\; \alpha
m_0+i\;B D {\bf m}'\times A {\bf m}\;,\nonumber\\
\beta B+\alpha D=B(1+t A)\;,\;\beta B+t A=B (1 + \alpha D),\;
B^2+AD=B(1+AD)\; ; \nonumber\end{eqnarray}
\begin{eqnarray}{\bf n}''=A\; {\bf m}''\qquad \Longrightarrow\qquad\qquad\qquad\nonumber\\
\beta m'_0\;A {\bf m}+B{\bf m}'\; \alpha m_0+i\; B {\bf m}'\times
A {\bf m}+\alpha m_0'\; {\bf m}+A {\bf
m}'\; m_0+i\; A {\bf m}'\times {\bf m}\nonumber\\
=A m'_0\; {\bf m}+A {\bf m}'\; m_0+i\; A {\bf m}' \times {\bf m}+A
t m_0'\; A {\bf m}+A D {\bf m}'\; \alpha
m_0+i\;A D {\bf m}'\times A {\bf m}\;,\nonumber\\
\beta A+\alpha=A(1+t A)\;,\; B\alpha+A=A(1 + \alpha D )\;,\; A(1+
B)=A(1+AD)\; ; \nonumber\end{eqnarray}
\begin{eqnarray}{\bf l}''=D\; {\bf m}''\qquad \Longrightarrow\qquad\qquad\qquad\nonumber\\
t m'_0\;B {\bf m}+D {\bf m}'\;\beta m_0+i\; D {\bf m}'\times B
{\bf m}+m_0'\;D {\bf m}+{\bf m}'\; t m_0 +
i\; {\bf m}'\times D {\bf m }\nonumber\\
= D m'_0\; {\bf m}+D {\bf m}'\; m_0+i\; D {\bf m}' \times {\bf
m}+D t m_0'\; A {\bf m}+D D {\bf m}'\; \alpha
m_0+i\; D D {\bf m}'\times A {\bf m}\;,\nonumber\\
tB+D=D(1+t A)\;,\; D\beta+t=D(1+\alpha D)\; ,\; D(1+ B)=D(1+AD)\;.
\nonumber\end{eqnarray}
By collecting together these results, we obtain the equations
\begin{eqnarray}\alpha (\beta+1)=\alpha (1+\alpha t)\;,\qquad A(B+1)
\alpha(1 +AD)\;,\nonumber\\
=\beta^2+\alpha t=\beta (1+\alpha t)\;,\qquad B^2 +AD
=\beta (1+AD)\;,\nonumber\\
t (\beta +1)=t(1+\alpha t)\;,\qquad D(B +1)=t(1+AD)
\;,\nonumber\\
\beta B+\alpha D=B(1+t A)\;,\;\beta B+t A=B (1
+ \alpha D),\; B^2+AD=B(1+AD)\;,\nonumber\\
\beta A+\alpha=A(1+t A)\;,\; B\alpha+A=A(1 +
\alpha D )\;,\; A(1+ B)=A(1+AD)\;,\nonumber\\
tB+D=D(1+t A)\;,\; D\beta+t=D(1+\alpha D)\; ,\; D(1+
B)=D(1+AD)\;.\nonumber\\\label{B4.3}\end{eqnarray}

The system (\ref{B4.3}) has a trivial solution with three
vanishing blocks:

\vspace{3mm}

\underline{solution $(M-1)$},
\begin{eqnarray}A=\alpha =0\;,\; B=\beta =0\;,\; D=t =0\;,\qquad G=\left(
    \ba{cc}0&0\\0&M\ea\right).\label{B4.4}\end{eqnarray}

Let us construct the solutions with two vanishing blocks. There
are three possible variants. The first is
\begin{eqnarray}A= \alpha=0\;,\qquad D=t=0\; ;\label{B4.5ax}\end{eqnarray}
the system (\ref{B4.3}) gives
\begin{eqnarray}0=0\;,\; 0=0\;,\; \beta^2=\beta\;,\; B^2=\beta\;,\;
 0=0\;,\;  0=0\;,\;
\beta B =B\;,\nonumber\\
 \beta B=B\;,\; B^2=B\;,  0=0\;,\; 0= 0\;,\; 0=0\;,\;
 0= 0\;,\; 0=0\;,\; 0=0\; ,
 \nonumber
 \label{B4.5bx}\end{eqnarray}
that is
\begin{eqnarray}\beta^2=\beta\;,\qquad B^2=\beta\;,\qquad
\beta B= B\;,\qquad B^2=B\;.\label{B4.5cx}\end{eqnarray}
The equations (\ref{B4.5cx}) have only already known solutions:
\vspace{2mm}

\underline{solution $(M-1)$},
\begin{eqnarray}A=\alpha =0\;,\; B=\beta =0\;,\; D=t =0\;,\qquad G=\left(
 \ba{cc}0&0\\0&M\ea\right) ;\label{B4.6}\end{eqnarray}

\underline{solution $(M-2)$},
\begin{eqnarray}A= \alpha=0\;,\; B=\beta=1\;,\; D=t=0\;,\qquad
 G=\left(\ba{cc}M&0\\0&M\ea\right).\label{B4.7}\end{eqnarray}

Now let us examine the variants with two vanishing blocks. The
first is
\begin{eqnarray}A= \alpha=0\;,\qquad B=\beta=0\; ;\label{B4.8ax}\end{eqnarray}
the system (\ref{B4.3}) gives
\begin{eqnarray}t=0\;,\; D= 0\;,\;
 0=0\;,\;0=0\;,\; 0=0\;,\; D=t\;,\;
 0=0\;,\; 0=0\;,\; 0=0\;,\nonumber\\
 0= 0\;,\; 0=0\;,\; 0=0\;,\; 0=0\;,\; t =D\;,\; 0=D\; ,
 \nonumber
 \label{B4.8bx}\end{eqnarray}
so we get only one

\vspace{3mm}

\underline{solution $(M-3)$},
\begin{eqnarray}A=\alpha =0\;,\; B=\beta =0\;,\; D=t\;,\nonumber\\
 G=\left(\ba{cc}0&0\\D M&M\ea\right),\;
 G'G=\left(\ba{cc}0&0\\DM' M&M'M\ea\right).
 \label{B4.8cx}\end{eqnarray}
This is a set of degenerate matrices of rank 2 with semi-group
structure.

The second set is
\begin{eqnarray}D= t=0\;,\qquad B=\beta=0\; ;\label{B4.9a}\end{eqnarray}
the system (\ref{B4.3}) reads
\begin{eqnarray}0=0\;,\; A=\alpha\;,\;
 0=0\;,\; 0 =0\;,\; 0=0\;,\;0=0\;,\nonumber\\
 0=0\;,\; 0=0,\; 0=0\;,\; \alpha=A\;,\; 0=0\;,
\; 0=0\;,\; 0=0\;,\; 0=0\;,\; 0= 0\;, \nonumber
\label{B4.9b}\end{eqnarray}
so we get

\vspace{3mm}

\underline{solution $(M-4)$},
\begin{eqnarray}A= \alpha\;,\; B=\beta=0\;,\; D=t=0\;,\qquad
 G=\left(\ba{cc}0&AM\\0&M\ea\right).\label{B4.9c}\end{eqnarray}

Now let us examine the solutions with one vanishing block. There
arise 3 possibilities.
    The first is
\begin{eqnarray}A=\alpha =0\; ;\label{B4.10ax}\end{eqnarray}
the system (\ref{B4.3}) gives
\begin{eqnarray}0=0\;,\; 0= 0\;,\;
\beta^2=\beta\;,\; B^2=\beta\;,\; t\beta=0\;,\; D(B +1)=t\;,\;
\beta B =B\;,\;\beta B=B \; , \nonumber\\
 B^2=B\;,\;
0= 0\;,\; 0=0\;,\; 0=0\;,\; tB =0\;,\; D\beta+t=D\;,\;D B=0\;.
\nonumber \label{B4.10b}\end{eqnarray}
This system has two solutions:

\vspace{3mm}

\underline{solution $(M-2)$},
\begin{eqnarray}A=\alpha =0\;,\; B=\beta=+1,\; D=t=0\;,\qquad
 G=\left(\ba{cc}M&0\\0&M\ea\right) ;\label{B4.10c}\end{eqnarray}

\vspace{3mm} \underline{solution $(M-3)$},
\begin{eqnarray}A=\alpha =0\;,\; B=\beta =0\;,\; D=t\;,\nonumber\\
G=\left(\ba{cc}0&0\\D M&M\ea\right),\; G'G=\left(\ba{cc} 0&0\\DM'
M&M'M\ea\right).\label{B4.10d}\end{eqnarray}

The third possibility is
\begin{eqnarray}B=\beta =0\;;\label{B4.11a}\end{eqnarray}
the system (\ref{B4.3}) reads
\begin{eqnarray}
0 =\alpha^2 t\;,\; A=\alpha\;,\;
\alpha t=0\;,\; AD= 0\;,\nonumber\\
0 =\alpha t^2\;,\; D=t\;,\;
 \alpha D=0\;,\;t A=0,\; AD=0\;,\nonumber\\
0=t A^2\;,\; 0=A^2 D\;,\; 0=A^2D\;,\;
 0=tD A\;,\; 0=\alpha D^2\;,\; 0=AD^2\;.
 \nonumber
 \label{B4.11b}\end{eqnarray}
These equations (\ref{B4.11b}) have two solutions:

\vspace{3mm}

\underline{solution $(M-3)$},
\begin{eqnarray}A=\alpha =0\;,\; B=\beta =0\;,\; D=t\;,\nonumber\\
 G=\left(\ba{cc}0&0\\D M&M\ea\right),\; G'G=\left(\ba{cc}
0&0\\DM' M&M'M\ea\right);\label{B4.11c}\end{eqnarray}

\underline{solution $(M-4)$},
\begin{eqnarray}A= \alpha\;,\; B=\beta=0\;,\; D=t=0\;,\;
 G=\left(\ba{cc}0&AM\\0&M\ea\right).\label{B4.11d}\end{eqnarray}

Consider the last variant with one vanishing block
\begin{eqnarray}D= t =0\; ;\label{B4.12a}\end{eqnarray}
the system (\ref{B4.3}) gives
\begin{eqnarray}\alpha\beta=0\;,\; A(B+1)=\alpha\;,\;
\beta^2=\beta\;,\; B^2=\beta\;,\nonumber\\
0=0\;,\; 0=0\;,\; \beta B =B\;,\;\beta B=B,\; B^2=B\;,\nonumber\\
\beta A+\alpha=A\;,\; B\alpha=0\;,\; A B=0\;,\; 0=0\;,\; 0 =0\;,\;
0=0\;.\label{B4.12b}\end{eqnarray}
These equations have two solutions:

\vspace{3mm}

\underline{solution $(M-2)$},
\begin{eqnarray}A= \alpha=0\;,\; B=\beta=1\;,\; D=t=0\;,\;
 G=\left(\ba{cc}M&0\\0&M\ea\right);\label{B4.12c}\end{eqnarray}

\underline{solution $(M-4)$},
\begin{eqnarray}A= \alpha\;,\; B=\beta=0\;,\; D=t=0\;,\;
 G=\left(\ba{cc}0&AM\\0&M\ea\right).\label{B4.12d}\end{eqnarray}

Now, let us consider the case where all the blocks are
non-vanishing
\begin{eqnarray}A\;,\; \alpha \neq 0\;,\qquad B\;,\;\beta \neq 0\;,
\qquad D\;,\; t \neq 0\; ;\label{B4.13a}\end{eqnarray}
The system (\ref{B4.3}) reads
\begin{eqnarray}\beta=\alpha t\;,\qquad A(B+1)=\alpha(1 +AD)\;,\nonumber\\
\beta^2+\alpha t=\beta (1+\alpha t)\;,\qquad B^2 +AD=\beta (1+AD)\;,\nonumber\\
\beta=\alpha t\;,\qquad D(B +1)=t(1 +AD)\;,\nonumber\\
\beta B+\alpha D=B(1+t A)\;,\;\beta B+t A=B (1+ \alpha D),\; B^2+AD=B(1+AD)\;,\nonumber\\
\beta A+\alpha=A(1+t A)\;,\qquad B=A D )\;,\qquad B=AD\;,\nonumber\\
B= AD\;,\qquad D\beta+t=D(1+\alpha D)\;,\qquad
B=AD\;,\nonumber\label{B4.13b}\end{eqnarray}
which is equivalent to
\begin{eqnarray}\beta=\alpha t\;,\qquad B=AD\;,\;
(D- t) (1+AD)=0\;,
\nonumber\\
  (A-\alpha) (1+AD)=0\;,\; (\alpha t-AD) (1+AD)=0\;,\nonumber\\
 (A-\alpha)(1+t A)=0\;,\qquad (D-t) (1+\alpha D)=0\;.\nonumber
 \label{B4.13c}\end{eqnarray}

There exist two possibilities. The first is
\begin{eqnarray}\beta=\alpha t\;,\; B=AD\;,\;
    D =-{1 \over A}\;,\; (A-\alpha)(1+t A)=0\; ;\label{B4.14a}\end{eqnarray}
we get two solutions:\par\medskip
\underline{solution $(M-5)$},
\begin{eqnarray}A=\alpha\;,\qquad B=-1,\qquad D=-{1 \over A}\;,\qquad\beta=At\;,\nonumber\\
G=\left(\ba{cc}(A t m_0-{\bf m}\vec{\sigma})&(A m_0+A {\bf m}\vec{\sigma})\\
(t m_0-A^{-1} {\bf m}\vec{\sigma})&(m_0+{\bf m}\vec{\sigma})
\ea\right);\label{B4.14b}\end{eqnarray}

\underline{solution $(M-6)$},
\begin{eqnarray}A\;, \alpha \neq 0\;,\qquad D=t=-{1 \over A}\;,\qquad
B=-1\;,\;\beta=-{\alpha \over A}\,,\nonumber\\G=\left(\ba{cc}
(-\alpha A^{-1} t m_0-{\bf m}\vec{\sigma})&(\alpha m_0+A {\bf m}\vec{\sigma})\\
(A^{-1} m_0+A^{-1} {\bf m}\vec{\sigma})&(m_0 +{\bf m}\vec{\sigma})
\ea\right).\label{B4.14c}\end{eqnarray}
The second possibility is
\begin{eqnarray}\beta=\alpha t\;,\qquad B=AD\;,\qquad D=t\;,\qquad A= \alpha\;;\label{B4.15a}\end{eqnarray}
it leads to

\vspace{3mm}

\underline{solution $(M-7)$},
\begin{eqnarray}A=\alpha \neq 0,\qquad D=t \neq 0,\qquad B=\beta=AD,\nonumber\\
G=\left(\ba{cc} AD M&A M\\D M&M\ea\right),\; G'G=\left(\ba{cc} AD
M'&A M'\\D M'&M'\ea\right)
\left(\ba{cc} AD M&A M\\D M&M\ea\right)\nonumber\\
=\left(\ba{cc} AD (1+ DA) M'M&A (1+DA) M'M\\D (1+DA) M'M&(1+ DA)
M'M\ea\right). \label{B4.15b}\end{eqnarray}
\section{One Independent Vector: Variant I(n)}
Consider {\bf variant $I(n)$}:
\begin{eqnarray}{\bf k}=A\; {\bf n}\;,\; k_0=\alpha\; n_0\;,\;
 {\bf m}=B\; {\bf n}\;,\; m_0=\beta\; n_0\;,\; {\bf l}=D\;
 {\bf n}\;,\; l_0=t\; n_0\;.\label{B5.1x}\end{eqnarray}
The multiplication law gives
\begin{eqnarray}\left.\ba{l} k_0''=\alpha^2\; n_0' n_0 +
A^2 {\bf n}' \cdot{\bf n} +t\; n'_0 n_0+D\; {\bf n}' \cdot{\bf n}\;,\\[2mm]
m_0''=\beta^2\; n_0' n_0+B^2 {\bf n}'\cdot {\bf n}
+t\; n'_0 n_0+D\; {\bf n}' \cdot{\bf n}\;,\\[2mm]
n_0''=\alpha\; n_0' n_0+A\; {\bf n}'\cdot {\bf n}
+\beta\; n'_0 n_0+B\; {\bf n}' \cdot{\bf n}\;,\\[2mm]
l_0''=t\alpha\; n_0' n_0+D A\; {\bf n}' \cdot{\bf n} +\beta t\;
n'_0 n_0+B D\; {\bf n}' \cdot{\bf n}\;,
 \ea \right.\nonumber\end{eqnarray}
\begin{eqnarray}\left.\ba{l} {\bf k}''=\alpha A\; n'_0 {\bf n} +
\alpha A\; {\bf n}' n_0 +i A^2\; {\bf n}'\times {\bf n} +
D\; n_0' {\bf n}+t\; {\bf n}' n_0+i D\; {\bf n}' \times {\bf n}\;,\\[2mm]
{\bf m}''=\beta B\; n'_0 {\bf n}+\beta B\; {\bf n}'n_0 +i B^2\;
{\bf n}'\times {\bf n} +
t\; n_0' {\bf n}+D\; {\bf n}' n_0+i D\; {\bf n}'\times {\bf n}\;,\\[2mm]
{\bf n}''=\alpha\; n'_0 {\bf n}+A\; {\bf n}' n_0+iA \; {\bf
n}'\times {\bf n}+B\; n_0' {\bf n}+\beta\; {\bf
n}' n_0+iB\; {\bf n}'\times {\bf n}\;,\\[2mm]
{\bf l}''=t A\; n'_0 {\bf n}+D \alpha\; {\bf n}' n_0+ i D A\;{\bf
n}'\times {\bf n} +\beta D\; n_0' {\bf n}+B t\;{\bf n}' n_0+i BD\;
{\bf n}'\times {\bf n }\;.\ea \right.\nonumber\end{eqnarray}
We further require

\vspace{3mm}

$\underline{\vspace{5mm} k_0''=\alpha\; n''_0: } $
\begin{eqnarray}(\alpha^2+t)\; n_0' n_0+(A^2+D)\; {\bf n}' \cdot{\bf n}
=\alpha (\alpha+\beta)\; n_0' n_0+\alpha (A+ B)\; {\bf n}' \cdot{\bf n}\;,\nonumber\\
(\alpha^2+t)=\alpha (\alpha+\beta)\;,\qquad (A^2+D)=\alpha (A+
B)\; ;\label{B5.2a}\end{eqnarray}

$\underline{\vspace{5mm} m_0''=\beta\; n''_0: } $
\begin{eqnarray}(\beta^2+t)\; n_0' n_0+(B^2 +D)\; {\bf n}'\cdot {\bf n}
=\beta\; (\alpha+\beta)\; n_0' n_0+\beta (A+B)\; {\bf n}'\cdot {\bf n}\;,\nonumber\\
(\beta^2+t)=\beta\; (\alpha+\beta)\;,\qquad (B^2 +D) =\beta
(A+B)\; ;\label{B5.2b}\end{eqnarray}

$\underline{\vspace{5mm} l_0''=t\; n''_0: } $
\begin{eqnarray}(t\alpha+\beta t)\; n_0' n_0+(D A+BD)\; {\bf n}'\cdot
{\bf n}=t\; (\alpha+\beta)\; n_0' n_0+t(A+B)\;{\bf n}'\cdot {\bf n}\;,\nonumber\\
(t\alpha+\beta t)=t\; (\alpha+\beta)\;,\qquad (D A +BD)=t(A+B)\;
;\label{B5.2c}\end{eqnarray}

$\underline{{\bf k} ''= A\; {\bf n}'' :} $
\begin{eqnarray}(\alpha A+D)\; n'_0 {\bf n}+(\alpha A +t)\; {\bf n}' n_0
+i (A^2+D)\; {\bf n}'\times {\bf n}\nonumber\\
= A(\alpha+B)\; n'_0 {\bf n} + A (A+\beta)\; {\bf n}' n_0 +
 i A(A +B)\; {\bf n}'\times {\bf n}\;,\nonumber\\
(\alpha A+D)= A(\alpha+B),\; (\alpha A +t)=A (A + \beta)\;,\;
(A^2+D)=A(A +B)\; ;\nonumber\\\label{B5.3a}\end{eqnarray}

$\underline{{\bf m} ''= B\; {\bf n}'' : }$
\begin{eqnarray}(\beta B+t)\; n'_0 {\bf n}+(\beta B+D)\; {\bf n}' n_0
+i (B^2+D)\; {\bf n}'\times {\bf n}\nonumber\\
=B(\alpha+B)\; n'_0 {\bf n} + B (A+\beta)\; {\bf n}' n_0 +
 i B(A +B)\; {\bf n}'\times {\bf n}\;,
\nonumber\end{eqnarray}
\begin{eqnarray}(\beta B+t)= B(\alpha+B)\;,\;(\beta B+D)= B (A +
\beta)\;,\; (B^2+D)=B(A +B)\;
;\nonumber\\\label{B5.3b}\end{eqnarray}

$\underline{{\bf l} ''= D\; {\bf n}'' :}$
\begin{eqnarray}(t A+\beta D)\; n'_0 {\bf n}+(D \alpha+B t)\; {\bf
n}' n_0+ i (D A+BD)\;{\bf n}'\times {\bf n}\nonumber\\
= D(\alpha+B)\; n'_0 {\bf n}+D (A+\beta)\; {\bf n}' n_0 + i D(A
+B)\; {\bf n}'\times {\bf n}\;,\nonumber\end{eqnarray}
\begin{eqnarray}(t A+\beta D)=D(\alpha+B)\;,\; (D \alpha+B t)
= D (A+\beta)\;,\;(D A+BD)=D(A
+B)\;.\nonumber\\\label{B5.3c}\end{eqnarray}

Thus, we have obtained the system
\begin{eqnarray}(\alpha^2+t)=\alpha (\alpha+\beta)\;,\qquad (A^2+D)=\alpha (A+ B)\;,\nonumber\\
(\beta^2+t)=\beta\; (\alpha+\beta)\;,\qquad (B^2 +D)=\beta (A+B)\;,\nonumber\\
t(\alpha+\beta)=t\; (\alpha+\beta)\;,\qquad D (A+B)= t(A+B)\;,\nonumber\\
(\alpha A+D)= A(\alpha+B)\;,\; (\alpha A +t)=A (A +\beta)\;,\; (A^2+D)=A(A +B)\;,\nonumber\\
(\beta B+t)= B(\alpha+B)\;,\; (\beta B+D)= B (A +\beta)\;,\; (B^2+D)=B(A +B)\;,\nonumber\\
(t A+\beta D)=D(\alpha+B)\;,\; (D \alpha+B t)= D (A+\beta)\;,\; D
(A+B)=D(A +B)\;, \nonumber\end{eqnarray}
which is equivalent to
\begin{eqnarray}
t=\alpha\beta\;,\; (A^2+D)=\alpha (A+ B)\;,\;
(B^2 +D)=\beta (A+B)\;,\nonumber\\
(D-t) (A+B)=0\;,\;  D=A B,\; (\alpha A +t)=A (A+\beta)\;,\nonumber\\
(\beta B+t)= B(\alpha+B)\;,\; (t A+\beta D)=D(\alpha+B)\;,\; (D
\alpha+B t)= D (A+\beta)\;. \nonumber \label{B5.4b}\end{eqnarray}
By eliminating the variables $D$ and $t$, we get
\begin{eqnarray}D= AB\;,\qquad t=\alpha\beta\; ;\label{B5.5a}\end{eqnarray}
the remaining independent equations are
\begin{eqnarray}(A-\alpha) (A+B)=0\;,\nonumber\\
(B-\beta) (B +A)=0\;,\nonumber\\ (AB- \alpha\beta) (A+B)=0\;,\nonumber\\
(A- \alpha) (A+\beta)=0\;,\nonumber\\ (B-\beta) (B+\alpha)= 0\;,\nonumber\\
A (\beta -B) (\alpha+B)=0\;,\nonumber\\  B (\alpha-A)
(A+\beta)=0\;.\label{B5.5b} \end{eqnarray}
First, assume that
\begin{eqnarray}A+B \neq 0\;,\label{B5.6a}\end{eqnarray}
then (\ref{B5.5b}) give
\begin{eqnarray}A=\alpha\;,\qquad B=\beta\;,\qquad  B (A-B) =0\;.\label{B5.6b} \end{eqnarray}
Here, there arise two solutions:

\vspace{3mm}

\underline{solution $(N-1)$},
\begin{eqnarray}A=\alpha\;,\; B=\beta=0,\; D=t =0,\qquad G=
    \left(\ba{cc}AN&N\\0&0\ea\right) ;\label{B5.7}\end{eqnarray}

\underline{solution $(N-2)$},
\begin{eqnarray}A=\alpha\;,\; B=\beta=A,\; D=t =A^2,\qquad G=\left(\ba{cc}
AN&N\\A^2 N&A N\ea\right),\nonumber\\
G'G=\left(\ba{cc}AN'&N'\\
A^2 N'&A N'\ea\right)\left(\ba{cc}A N&N\\A^2 N&A N\ea\right)=
\left(\ba{cc}A\; 2A N'N&2AN'N\\A^2\; 2AN'N&A\; 2A N'N
\ea\right).\label{B5.8}\end{eqnarray}

Now, let it be
\begin{eqnarray}A+ B=0\;,\qquad B=-A\; ;\label{B5.9a}\end{eqnarray}
correspondingly to (\ref{B5.5b}) it follows
\begin{eqnarray} B =-A\;,\qquad D= -A^2\;,\qquad t=\alpha\beta\;,\nonumber\\
(A- \alpha) (A+\beta)=0\;,\nonumber\\ A (\beta +A) (\alpha-A)=0\;,\nonumber\\
 A (\alpha-A)  (A+\beta)=0\;.\label{B5.9b} \end{eqnarray}
Two solutions are possible:

\underline{solution $(N-3)$},
\begin{eqnarray}\beta =B =-A\;,\qquad D= -A^2\;,\qquad t=-\alpha A\;,\nonumber\\
G=\left(\ba{cc}\alpha n_0+A{\bf n}\vec{\sigma}&n_0+{\bf n}\vec{\sigma}\\[2mm]
-\alpha A n_0-A^2 {\bf n}\vec{\sigma}&-A n_0-A{\bf n}\vec{\sigma}
\ea\right),\qquad\det\; G=0\; ;\label{B5.10}\end{eqnarray}

\underline{solution $(N-4)$},
\begin{eqnarray}A =+\alpha\;,\qquad B =-A\;,\qquad D= -A^2\;,\qquad t = A\beta\;,\nonumber\\
G=\left(\ba{cc} A n_0+A{\bf n}\vec{\sigma}&n_0+{\bf n}\vec{\sigma}\\[2mm]
\beta A n_0-A^2 {\bf n}\vec{\sigma}&\beta n_0-A{\bf n}\vec{\sigma}
\ea\right),\qquad\det\; G=0\;.\label{B5.11}\end{eqnarray}
\section{One Independent Vector: Variant I(l)}
Let us examine the {\bf variant $I(l)$}:
\begin{eqnarray}{\bf k}=A\; {\bf l}\;,\; k_0=\alpha\; l_0\;,\;
{\bf m}=B\; {\bf l}\;,\; m_0=\beta\; l_0\;,\; {\bf n}=D\; {\bf
l}\;,\; n_0=t\; l_0\;.\label{B6.1}\end{eqnarray}
The multiplication law takes the form
\begin{eqnarray}k_0''=\alpha^2 l_0'\; l_0+A^2 {\bf l}'\cdot {\bf l}
+ t l'_0\; l_0+D {\bf l}'\cdot{\bf l}\;,\nonumber\\
m_0''=\beta^2 l_0'\; l_0+B^2 {\bf l}'\cdot {\bf l} +t l'_0\; l_0 +
 D {\bf l}'\cdot {\bf l}\;,\nonumber\\
n_0''=\alpha t l_0'\; l_0+A D {\bf l}'\cdot {\bf l}
+\beta t l'_0\; l_0+B D {\bf l}'\cdot {\bf l}\;,\nonumber\\
l_0''=\alpha l_0'\; l_0+A {\bf l }'\cdot{\bf l} +\beta l'_0\;
l_0+B {\bf l}'\cdot {\bf l}\;,\label{B6.2a}\end{eqnarray}
\begin{eqnarray}{\bf k}''=A \alpha l'_0\; {\bf l} +
 A\alpha {\bf l}'\; l_0+i\; A^2 {\bf l}'\times {\bf l}+t l_0'\; {\bf l}+D {\bf\; l}'\; l_0 +
i\;D {\bf l}'\times {\bf l}\;,\nonumber\\
{\bf m}''=B\beta l'_0\; {\bf l}+B\beta {\bf l}'
\;l_0+i\; B^2 {\bf l}'\times {\bf l}+D l_0'\; {\bf l}+t {\bf l}'\; l_0+i\;D {\bf l}'\times {\bf l}\;,\nonumber\\
{\bf n}''=\alpha D l'_0\; {\bf l}+A t {\bf l}'\; l_0 + i\;A D {\bf
l}'\times {\bf l}+B t l_0'\; {\bf l}+\beta D
{\bf l}'\; l_0+i\;BD {\bf l}'\times {\bf l}\;,\nonumber\\
{\bf l}''=A l'_0\; {\bf l}+\alpha {\bf l}'\; l_0+i \;A {\bf
l}'\times {\bf l}+\beta l_0'\; {\bf l}+B {\bf l}'\; l_0+i\;B {\bf
l}'\times {\bf l }\;.\nonumber\\\label{6.2bz}\end{eqnarray}

We further require
\begin{eqnarray}k_0''=\alpha\; l_0 ''\qquad \Longrightarrow\nonumber\\
(\alpha^2+t) l_0'\; l_0+(A^2+D) {\bf l}'\cdot {\bf l}=\alpha
(\alpha+\beta) l_0'\; l_0+\alpha ( A+B) {\bf l}'\cdot {\bf l}
\;\; \Longrightarrow\nonumber\\
(\alpha^2+t)=\alpha (\alpha+\beta)\;,\qquad (A^2 + D)=\alpha
(A+B)\; ; \nonumber\end{eqnarray}
\begin{eqnarray}m_0''=\beta\; l_0''\qquad \Longrightarrow\nonumber\\
(\beta^2+t) l_0'\; l_0+(B^2+D) {\bf l}'\cdot {\bf l}=\beta
(\alpha+\beta) l_0'\; l_0+\beta (A+B) {\bf l}'\cdot {\bf l}
 \;\; \Longrightarrow\nonumber\\
(\beta^2+t)=\beta (\alpha+\beta)\;,\qquad (B^2 + D)=\beta (A+B)\;
; \nonumber\end{eqnarray}
\begin{eqnarray}n_0''=t\; l_0''\qquad \Longrightarrow\nonumber\\
(\alpha t+\beta t) l_0'\; l_0+(A D+B D) {\bf l}' \cdot {\bf l} =t
(\alpha+\beta) l_0'\; l_0+t (A+B) {\bf l}'\cdot {\bf l}
\;\; \Longrightarrow\nonumber\\
(\alpha t+\beta t)= t (\alpha+\beta)\;,\qquad (A D + B D)=t
(A+B)\; ; \nonumber\end{eqnarray}
\begin{eqnarray}{\bf k} ''=A\; {\bf l}''\qquad \Longrightarrow\nonumber\\
(A \alpha+t) l'_0\; {\bf l} +
 (A\alpha+D) {\bf l}'\; l_0+i\; (A^2 +D) {\bf l}'\times {\bf l}\nonumber\\
=
 A (A+\beta) l'_0\; {\bf l}+A (\alpha+B) {\bf l}'\;
l_0+i\; A (A +B) {\bf l}'\times {\bf l}\;
\Longrightarrow\nonumber\\
(A \alpha+t)=A (A+\beta)\;,\qquad (A\alpha+D)=A
(\alpha+B)\;,\qquad (A^2 +D)=A (A +B)\; ; \nonumber\end{eqnarray}
\begin{eqnarray}{\bf m}''=B\; {\bf l}''\qquad \Longrightarrow\nonumber\\
(B\beta +D) l'_0\; {\bf l}+(B\beta+t) {\bf l}'
\;l_0+i\; (B^2+D){\bf l}'\times {\bf l}\nonumber\\
= B (A+\beta) l'_0\; {\bf l}+B (\alpha +B) {\bf l}'\; l_0+i\;B
(A+B) {\bf l}'\times {\bf l}\;\;
\Longrightarrow\nonumber\\
(B\beta +D)=B (A+\beta)\;,\qquad (B\beta+t)=B ( \alpha
+B)\;,\qquad (B^2+D)=B (A+B)\; ; \nonumber\end{eqnarray}
\begin{eqnarray}{\bf n}''=D\; {\bf l} ''\qquad \Longrightarrow\nonumber\\
(\alpha D +B t) l'_0\; {\bf l}+(A t+\beta D ) {\bf l}'
\; l_0+i\; (A D+BD) {\bf l}'\times {\bf l}\nonumber\\
=
 D (A+\beta) l'_0\; {\bf l} +
 D (\alpha+B) {\bf l}'\; l_0+i\;D (A +B) {\bf l}'\times {\bf l}\;\; \Longrightarrow\nonumber\\
(\alpha D +B t)=D (A+\beta)\;,\; (A t+\beta D) = D (\alpha+B)\;,\;
(A D+BD)=D (A +B)\;. \nonumber\end{eqnarray}
By collecting the results together, we get
\begin{eqnarray}(\alpha^2+t)=\alpha (\alpha+\beta)\;,\qquad (A^2
+ D)=\alpha (A+B)\;,\nonumber\\
(\beta^2+t)=\beta (\alpha+\beta)\;,\qquad (B^2
+ D)=\beta (A+B)\;,\nonumber\\
(\alpha t+\beta t)= t (\alpha+\beta)\;,\qquad (A D +
B D)=t (A+B)\;,\nonumber\\
(A \alpha+t)=A (A+\beta)\;,\; (A\alpha+D)=A
(\alpha+B)\;,\; (A^2 +D)=A (A +B)\;,\nonumber\\
(B\beta +D)=B (A+\beta)\;,\; (B\beta+t)=B (
\alpha +B)\;,\; (B^2+D)=B (A+B)\;,\nonumber\\
(\alpha D +B t)=D (A+\beta)\;,\; (A t+\beta D) = D (\alpha+B)\;,\;
(A D+BD)=D (A +B)\;.\nonumber\\\label{B6.3}\end{eqnarray}
The system (\ref{B6.3}) is equivalent to
\begin{eqnarray} t=\alpha\beta\;,\qquad
(A^2+D)=\alpha (A+B)\;,\nonumber\\
 t=\beta \alpha\;,\qquad
(B^2+D)=\beta (A+B)\;,\nonumber\\
0=0\;,\qquad (D- t)(A+B)=0\;,\nonumber\\
(A \alpha+t)=A (A+\beta)\;,\qquad (A\alpha+D)=A
(\alpha+B)\;,\qquad D=A B\;,\nonumber\\
(B\beta +D)=B (A+\beta)\;,\qquad (B\beta+t)=B ( \alpha
+B)\;,\qquad
 D=AB\;,\nonumber\\
(\alpha D +B t)=D (A+\beta)\;,\qquad (A t+\beta D) = D
(\alpha+B)\;,\qquad 0=0\;.\label{B6.4}\end{eqnarray}
By excluding the variables $t$ and $D$ via
\begin{eqnarray} t=\alpha\beta\;,\qquad D=AB\;,\label{B6.5a} \end{eqnarray}
we obtain
\begin{eqnarray} (A- \alpha)(A+B)=0\;,\; (B-\beta)(A+B)=0\;,\;
(AB- \alpha\beta)(A+B)=0\;, \nonumber\\ (A-\alpha) (A+\beta)
=0\;,\quad
 (B-\beta) (\alpha +B)=0\;,\qquad
 \nonumber\\ A(\beta-B) (\alpha+ B )=0\;,\quad
B (\alpha -A) (A +\beta)=0\;.\qquad \label{B6.5b}\end{eqnarray}
Because these equations coincide with (\ref{B5.5b}), we may use
yet the known results:

\vspace{3mm}

\underline{solution $(L-1)$},
\begin{eqnarray}A=\alpha\;,\; B=\beta=0\;,\; D=t =0\;,\qquad G
= \left(\ba{cc}AL&0\\L&0\ea\right),\nonumber\\
 G'G=\left(\ba{cc}AL&0\\L&0\ea\right)
\left(\ba{cc}AL'&L'\\0&0\ea\right)= \left(\ba{cc}
AL'L&0\\L'L&0\ea\right) ;\label{B6.6}\end{eqnarray}

\underline{solution $(L-2)$},
\begin{eqnarray}A=\alpha\;,\; B=\beta=A\;,\; D=t =A^2\;,\qquad G=\left(\ba{cc}
A L&A^2 L\\ L&A L\ea\right),\nonumber\\
G'G=\left(\ba{cc}AL'&A^2 L'\\L '&A L'\ea\right) \left(\ba{cc}A
L&A^2 L\\L&A L\ea\right)= \left(\ba{cc}A\; 2A L'L&A^2\; 2AL'L\\
2AL'L&A\; 2A L'L \ea\right) ;\label{B6.7}\end{eqnarray}

\underline{solution $(L-3)$},
\begin{eqnarray}\beta =B =-A\;,\qquad
 D= -A^2\;,\qquad t=-\alpha A\;,\nonumber\\
G=\left(\ba{cc}
\alpha l_0+A{\bf l}\vec{\sigma}&-\alpha A l_0-A^2 {\bf l}\\[2mm]
l_0+{\bf l}\vec{\sigma} &-A l_0-A{\bf l}\vec{\sigma}
\ea\right),\qquad\det\; G=0\; ;\label{B6.8}\end{eqnarray}

\underline{solution $(L-4)$},
\begin{eqnarray}A =+\alpha\;,\qquad B =-A\;,\qquad D= -A^2\;,\qquad t
= A\beta\;,\nonumber\\
G=\left(\ba{cc}
A l_0+A{\bf l}\vec{\sigma}&\beta Al_0-A^2 {\bf l}\vec{\sigma}\\[2mm]
l_0+{\bf l}\vec{\sigma} &\beta l_0-A{\bf l}\vec{\sigma}
\ea\right),\qquad\det\; G=0\;.\label{B6.9}\end{eqnarray}
\section{Two Independent Vectors: Variant II(k,m)}
Let us examine the {\bf variant II(k,m)}:
\begin{eqnarray}\left.\ba{ll} {\bf n}=A {\bf k}+B {\bf m}\;,
\qquad&n_0=\alpha k_0+\beta m_0\;,\\[2mm]
{\bf l}=C {\bf k}+D {\bf m}\;,\qquad&l_0=s k_0 +t m_0\;. \ea
\right.\label{B7.1}\end{eqnarray}
The multiplication law takes the form
\begin{eqnarray}\left.\ba{l} k_0''=k_0'\; k_0+{\bf k}'\cdot
{\bf k}+(\alpha k'_0+\beta m'_0) (s k_0+t m_0)+ (A {\bf k}'+B {\bf
m}'
)(C {\bf k}+D {\bf m})\;,\\[2mm]
m_0''=m_0'\; m_0+{\bf m}'\cdot {\bf m} +(s k'_0+t m'_0) (\alpha
k_0+\beta m_0 +
 (C {\bf k}'+D {\bf m}')(A {\bf k}+B {\bf m})\;,\\[2mm]
n_0''=k_0' (\alpha k_0+\beta m_0)+{\bf k}' (A {\bf k}+B {\bf m})
+(\alpha k'_0+\beta m'_0) m_0 +
 (A {\bf k}'+B {\bf m}') {\bf m}\;,\\[2mm]
l_0''=(s k'_0+t m'_0) k_0+(C {\bf k}'+D {\bf m}') {\bf k}
 +
 m'_0 (s k_0+t m_0) +
 {\bf m}' (C {\bf k}+D {\bf m})\;,
\ea \right. \nonumber\end{eqnarray}
\begin{eqnarray}{\bf k}''=k'_0\; {\bf k}+{\bf k}'\; k_0+i\; {\bf
k}'\times {\bf k}
 +(\alpha k'_0+\beta m'_0) (C {\bf k}
+ D {\bf m})\nonumber\\
+ (A{\bf k}'+B {\bf m}') (s k_0+t m_0)+i\; (A {\bf k}'
+ B {\bf m}')\times (C {\bf k}+D {\bf m})\;,\nonumber\\
{\bf m}''=m'_0\; {\bf m}+{\bf m}'\; m_0+i\; {\bf m}'\times {\bf m}
+ (s k'_0+t m'_0) (A {\bf k}+B {\bf m})\nonumber\\
+ (C {\bf k}'+D {\bf m}') (\alpha k_0+\beta m_0)
+i\; (C {\bf k}'+D {\bf m}')\times  (A {\bf k}+B {\bf m})\;,\nonumber\\
{\bf n}''=k'_0 (A {\bf k}+B {\bf m})+{\bf k}' (\alpha
k_0+\beta m_0) +i\; {\bf k}'\times (A {\bf k} +  B {\bf m}) \nonumber\\
 + (\alpha k'_0+\beta m'_0) {\bf m}+(A {\bf k}'+B {\bf m}') m_0 +
 i\;(A {\bf k}'+B {\bf m}')\times {\bf m}\;,\nonumber\\
{\bf l}''=(s k'_0+t m_0') {\bf k}+(C {\bf k}'+D
{\bf m}') k_0+i\; (C {\bf k}'+D {\bf m}')\times {\bf k}\nonumber\\
+ m_0' (C {\bf k}+D {\bf m})+{\bf m}' (s k_0+t m_0)+i\; {\bf
m}'\times (C {\bf k}+D {\bf m})\;. \nonumber\end{eqnarray}

By requiring $\alpha k''_0+\beta m_0 ''=n_0'' $, or
\begin{eqnarray}\alpha\; [\; k_0'\; k_0+{\bf k}'\cdot {\bf k}+(\alpha
k'_0+\beta m'_0) (s k_0+t m_0)+ (A {\bf k}'+B
{\bf m}')(C {\bf k}+D {\bf m})\; ]\;\nonumber\\
+\beta\; [\; m_0'\; m_0+{\bf m}'\cdot {\bf m} +(s k'_0+t m'_0)
(\alpha k_0+\beta m_0 +
 (C {\bf k}'+D {\bf m}')(A {\bf k}+B {\bf m})\; ]\nonumber\\
= k_0' (\alpha k_0+\beta m_0) + {\bf k}' (A {\bf k}+B {\bf m})
+(\alpha k'_0+\beta m'_0) m_0 +
 (A {\bf k}'+B {\bf m}') {\bf m}\;,\nonumber \end{eqnarray}
we obtain
\begin{eqnarray}\left.\ba{lll}
k'_0 k_0:\qquad&\alpha+\alpha^2 s+\alpha\beta s=\alpha\;,\quad
m'_0 m_0:\qquad&\alpha\beta t+\beta+t\beta^2=\beta\;,\\
k'_0 m_0:\qquad&\alpha^2 t+\beta^2 s=\beta+\alpha\;,\quad
m'_0 k_0:\qquad&\alpha\beta s+\alpha\beta t=0\;,\\
{\bf m}' \cdot{\bf k} :\qquad&\alpha BC+\beta AD=0\;,\quad
{\bf k}' \cdot{\bf m} :\qquad&\alpha AD+\beta CB=B+A\;,\\
{\bf k} '\cdot{\bf k} :\qquad&\alpha +\alpha AC +\beta
AC=A\;,\quad {\bf m}'\cdot {\bf m} :\qquad&\alpha BD+\beta+\beta
DB=B\;. \ea \right.\label{B7.2b}\end{eqnarray}

We further require $s k''_0+t m_0''=l''_0$, or
\begin{eqnarray}s\; [\; k_0'\; k_0+{\bf k}'\cdot {\bf k}+(\alpha k'_0 +
\beta m'_0) (s k_0+t m_0)+ (A {\bf k}'+B {\bf m}' )(C {\bf k}+D {\bf m})\; ]\;\nonumber\\
+ t\; [\; m_0'\; m_0+{\bf m}'\cdot {\bf m}+(s k'_0+t m'_0) (\alpha
k_0+\beta m_0 +
 (C {\bf k}'+D {\bf m}')(A {\bf k}+B {\bf m})\; ]\nonumber\\
=(s k'_0+t m'_0) k_0+(C {\bf k}'+D {\bf m}'){\bf k} + m'_0 (s
k_0+t m_0) +
 {\bf m}' (C {\bf k}+D {\bf m})\,,\nonumber \end{eqnarray}
and we find
\begin{eqnarray}\left.\ba{lll}
k'_0 k_0:\qquad&s+s^2 \alpha+t s \alpha=s\;,\quad
m'_0 m_0:\qquad&\beta s t+t+\beta t^2=t\;,\\
k'_0 m_0:\qquad&st \alpha+st\beta=0\;,\quad
m'_0 k_0:\qquad&\beta s^2+\alpha t^2=t+s\;,\\
{\bf k}' \cdot{\bf m} :\qquad&sAD+t CB=0,\quad
{\bf m}'\cdot {\bf k} :\qquad&s BC+t AD=C+ D\;,\\
{\bf k} '\cdot{\bf k} :\qquad&s +s AC +t AC=C\;,\quad {\bf m}'
\cdot {\bf m} :\qquad&s BD+t+t BD=D\;. \ea
\right.\label{B7.3b}\end{eqnarray}

By requiring $ A {\bf k}''+B {\bf m}''={\bf n}''$, or
\begin{eqnarray}A\; [\; k'_0\; {\bf k}+{\bf k}'\; k_0+i\; {\bf k}'\times {\bf k}
 + (\alpha k'_0+\beta m'_0) (C {\bf k}+D {\bf m})\nonumber\\+(A
{\bf k}'+B {\bf m}') (s k_0+t m_0)+i\; (A {\bf k}'
+ B {\bf m}')\times (C {\bf k}+D {\bf m})\; ]\nonumber\\
+ B\; [\; m'_0\; {\bf m}+{\bf m}'\; m_0+i\; {\bf m}'\times {\bf m}
+
 (s k'_0+t m'_0) (A {\bf k}+B {\bf m})\nonumber\\
 +(C {\bf k}'+D {\bf m}') (\alpha k_0+\beta m_0)+i\; (C {\bf k}'+D {\bf m}')\times
 (A {\bf k}+B {\bf m})\; ]\nonumber\\
= k'_0 (A {\bf k}+B {\bf m})+{\bf k}' (\alpha k_0 +
\beta m_0) +i\; {\bf k}'\times (A {\bf k} + B {\bf m})\nonumber\\
 + (\alpha k'_0+\beta m'_0) {\bf m}+(A {\bf k}'+B {\bf m}') m_0 +
 i\;(A {\bf k}'+B {\bf m}')\times {\bf m}\;,\nonumber\end{eqnarray}
we yield
\begin{eqnarray}k_0' {\bf k}:\qquad A+\alpha A C+s BA=A\;,\quad
{\bf m}' m_0:\qquad t AB+B+\beta BD=B\;,\nonumber\\
{\bf k}' k_0 :\qquad A+s A^2+\alpha BC=\alpha\;,\quad
m'_0 {\bf m}:\qquad\beta AD+B+tB^2=\beta\;,\nonumber\\
k'_0 {\bf m}:\qquad \alpha AD+sB^2=B+\alpha\;,\quad
{\bf k}' m_0:\qquad tA^2+\beta CB=\beta +A\;,\nonumber\\
m'_0 {\bf k}:\qquad\beta AC+t AB=0\;,\quad
+{\bf m}' k_0:\qquad s AB+\alpha BD=0\;,\nonumber\\
{\bf k}'\times {\bf k} :\qquad A+A^2 C+BCA=A\;,\quad
{\bf m}'\times {\bf m} :\qquad ABD+B+DB^2=B\;,\nonumber\\
{\bf k}'\times {\bf m} :\qquad A^2 D+CB^2=B +A\;,\quad {\bf
m}'\times {\bf k} :\qquad ABC+ADB=0\;. \nonumber
\label{B7.4b}\end{eqnarray}

Now, requiring $ C {\bf k}''+D {\bf m}''={\bf l}'' $, or
\begin{eqnarray}C\; [\; k'_0\; {\bf k}+{\bf k}'\; k_0+i\; {\bf k}' \times {\bf k}
 + (\alpha k'_0+\beta m'_0) (C {\bf k}+D {\bf m})\nonumber\\
+(A {\bf k}'+B {\bf m}') (s k_0+t m_0)+i\; (A {\bf k}'
+ B {\bf m}')\times (C {\bf k}+D {\bf m}) ]\nonumber\\
 +D\; [\; m'_0\; {\bf m}+{\bf m}'\; m_0+i\; {\bf m}'\times {\bf m}
+ (s k'_0+t m'_0) (A {\bf k}+B {\bf m})\nonumber\\
+(C {\bf k}'+D {\bf m}') (\alpha k_0+\beta m_0)
+i\; (C {\bf k}'+D {\bf m}')\times (A {\bf k}+B {\bf m})\nonumber\\
= (s k'_0+t m_0') {\bf k}+(C {\bf k}'+D {\bf m}')
k_0+i\; (C {\bf k}'+D {\bf m}')\times {\bf k}\nonumber\\
 + m_0' (C {\bf k}+D {\bf m})+{\bf m}' (s k_0+t
m_0)+i\; {\bf m}'\times (C {\bf k}+D {\bf
m})\;,\nonumber\end{eqnarray}
 we find
\begin{eqnarray}k_0' {\bf k}:\qquad C+\alpha C^2+s AD=s\;,\quad
{\bf m}' m_0:\qquad t CB+D+\beta D^2 =t\;,\nonumber\\
{\bf k}' k_0 :\qquad C+s AC+\alpha CD=C\;,\quad
m'_0 {\bf m}:\qquad\beta CD+D+t BD=D\;,\nonumber\\
k'_0 {\bf m}:\qquad \alpha DC+s BD =0\;,\quad
{\bf k}' m_0:\qquad t AC+\beta CD=0\;,\nonumber\\
m'_0 {\bf k}:\qquad\beta C^2+t AD=t+C\;,\quad
+{\bf m}' k_0:\qquad s CB+\alpha D^2=D +s\;,\nonumber\\
{\bf k}'\times {\bf k} :\qquad C +AC^2+DCA=C\;,\quad
{\bf m}'\times {\bf m} :\qquad CBD+D+B D^2=D\;,\nonumber\\
{\bf k}'\times {\bf m} :\qquad ACD+CBD =0\;,\quad {\bf m}'\times
{\bf k} :\qquad BC^2+AD^2=D+C\;. \nonumber
\label{B7.5b}\end{eqnarray}

We collect the results together:
\begin{eqnarray}\left.\ba{ll} \alpha s(\alpha+\beta)=0\;,
\qquad&\beta t (\alpha+\beta)=0\;,\\[2mm]
\alpha^2 t+\beta^2 s=\beta+\alpha\;,\qquad&\alpha\beta (s+ t)=0\;,\\[2mm]
\alpha BC+\beta AD=0\;,\qquad&\alpha AD+\beta BC=B+A\;,\\[2mm]
\alpha+(\alpha +\beta) AC=A\;,\qquad&\beta+(\alpha +\beta)
BD=B\;,\ea \right.\nonumber\end{eqnarray}
\begin{eqnarray} \left.\ba{ll}
 s \alpha (s +t) =0\;,\qquad&\beta t (s + t) =0\;,\\[2mm]
st (\alpha+\beta)=0\;,\qquad&\beta s^2+\alpha t^2=t+s\;,\\[2mm]
 sAD+t CB=0,\qquad&s BC+t AD=C+ D\;,\\[2mm]
s+AC(s+t)=C\;,\qquad&t+BD(t+s) =D\;,\ea
\right.\nonumber\end{eqnarray}
\begin{eqnarray} \left.\ba{ll}
A(\alpha C+s B)=0\;,\qquad&B (\beta D +t A)=0\;,\\[2mm]
 A+s A^2+\alpha BC=\alpha\;,\qquad&B+tB^2 +\beta AD=\beta\;,\\[2mm]
\alpha AD+sB^2=B+\alpha\;,\qquad&\beta BC+tA^2 =A+\beta\;,\\[2mm]
\beta AC+t AB=0\;,\qquad&\alpha BD+s AB =0\;, \ea \right.\nonumber\\
\left.\ba{ll}A+AC (A+B)=A\;,\qquad&B+BD (B+A) =B\;,\\[2mm]
 A^2 D+B^2 C=A+ B\;,\qquad&AB (C+D)=0\;,\\[2mm]
 C+\alpha C^2+s AD=s\;,\qquad&D+\beta D^2+t BC =t\;,\\[2mm]
C (s A+\alpha D)=0\;,\qquad&D(t B +\beta C )=0\;, \ea \right.\nonumber\\
\left.\ba{ll} D( \alpha C+s B) =0\;,\qquad&C(\beta D+t A)=0\;,\\[2mm]
\beta C^2+t AD=t+C\;,\qquad&\alpha D^2 +s BC=D +s\;,\\[2mm]
 AC(C+D)=0\;,\qquad&BD(D+C)=0\;,\\[2mm]
 DC ( A+B) =0\;, \qquad&BC^2+AD^2=D+C\;.
 \ea \right.\label{B7.6}\end{eqnarray}

While analyzing this system, we will assume that the blocks
\begin{eqnarray}(A,\alpha)\;,\qquad (B,\beta)\;,\qquad (D,t)\;,\qquad (C,s)\label{B7.7b}\end{eqnarray}
can be vanishing or not. The simplest possibility is
\begin{eqnarray}(A,\alpha)=0\;,\qquad (B,\beta)=0\;,\qquad (D,t)=0\;,
    \qquad (C,s)=0\;,\label{B7.8a}\end{eqnarray}
which leads to

\vspace{3mm}

 \underline{solution $(KM-1)$},
\begin{eqnarray}G=\left(\ba {cc} K&0\\ 0&M \ea\right)\;.\label{B7.8b}\end{eqnarray}

One may separate 4 different cases with one non-vanishing block:
\begin{eqnarray}(A,\alpha) \neq 0\;,\qquad (B,\beta)=0\;,\qquad (D,t)=0\;,\qquad (C,s)=0\; ;\nonumber\\
(A,\alpha)=0\;,\qquad (B,\beta) \neq 0\;,\qquad (D,t)=0\;,\qquad (C,s) = 0\; ;\nonumber\\
(A,\alpha)=0\;,\qquad (B,\beta) =0\;,\qquad (D,t) \neq 0\;,\qquad (C,s)=0\; ;\nonumber\\
(A,\alpha)=0\;,\qquad (B,\beta)=0\;,\qquad (D,t)=0\;,\qquad (C,s)
\neq 0\;.\label{B7.9a}\end{eqnarray}

Moreover, there are six possible  cases with two non-vanishing
blocks (two of them marked by
    asterisk were considered in previous sections)
\begin{eqnarray}(A,\alpha)=0\;,\qquad (B,\beta)=0\;,\qquad (D,t) \neq 0
\;,\qquad (C,s) \neq 0\; ;\nonumber\\
(A,\alpha)=0\;,\qquad (B,\beta) \neq 0\;,\qquad (D,t)=0
\;,\qquad (C,s) \neq 0\; ;\nonumber\\
(A,\alpha)=0\;,\qquad (B,\beta) \neq 0\;,\qquad (D,t) \neq
0\;,\qquad (C,s)=0\;;\qquad (*)\;\nonumber\\
(A,\alpha) \neq 0\;,\qquad (B,\beta)=0\;,\qquad (D,t) \neq
0\;,\qquad (C,s)=0\;;\qquad (*)\;\nonumber\\
(A,\alpha) \neq 0\;,\qquad (B,\beta)= 0\;,\qquad (D,t)=0\;,\qquad (C,s) \neq 0\; ;\nonumber\\
(A,\alpha) \neq 0\;,\qquad (B,\beta) \neq 0\;,\qquad (D,t) =
0\;,\qquad (C,s)=0\;.\label{B7.9b}\end{eqnarray}

There exist four  cases with three non-vanishing blocks:
\begin{eqnarray}(A,\alpha)=0\;,\qquad (B,\beta)\neq 0\;,\qquad (D,t) \neq 0
\;,\qquad (C,s) \neq 0\; ;\nonumber\\
(A,\alpha) \neq 0\;,\qquad (B,\beta)= 0\;,\qquad (D,t) \neq
0\;,\qquad (C,s) \neq 0\; ;\nonumber\\
(A,\alpha) \neq 0\;,\qquad (B,\beta)\neq 0\;,\qquad (D,t)=0\;,\qquad (C,s) \neq 0\; ;\nonumber\\
(A,\alpha) \neq 0\;,\qquad (B,\beta)\neq 0\;,\qquad (D,t) \neq
0\;,\qquad (C,s)=0\;.\label{B7.9c}\end{eqnarray}

Finally, there is possible the case when all  the four blocks are
non-zero.

Let it be (see (\ref{B7.9a}))
\begin{eqnarray}(A,\alpha) \neq 0\;,\qquad (B,\beta)=0\;,\qquad (D,t)=0
    \;,\qquad (C,s)=0\; ;\label{B7.10}\end{eqnarray}
the system (\ref{B7.6}) becomes simpler and leads to the already
known solution (\ref{B7.8b}).

Let it be (see (\ref{B7.9a}))
\begin{eqnarray}(A,\alpha)=0\;,\qquad (B,\beta) \neq 0\;,\qquad (D,t)=0
\;,\qquad (C,s)=0\; ;\label{B7.11}\end{eqnarray}
the system (\ref{B7.6}) again gives the known solution
(\ref{B7.8b}).

Let it be (see (\ref{B7.9a}))
\begin{eqnarray}(A,\alpha)=0\;,\qquad (B,\beta)= 0\;,\qquad (D,t) \neq
0\;,\qquad (C,s)=0\; ;\label{B7.12}\end{eqnarray}
again we arrive at the solution (\ref{B7.8b}).

Let it be (see (\ref{B7.9a}))
\begin{eqnarray}(A,\alpha)=0\;,\qquad (B,\beta)= 0\;,\qquad (D,t)=0
\;,\qquad (C,s) \neq 0\; ;\label{B7.13}\end{eqnarray}
it leads to the known solution (\ref{B7.8b}).

Now, let it be (see (\ref{B7.9b}))
\begin{eqnarray}(A,\alpha)=0\;,\qquad (B,\beta)=0\;,\qquad (D,t) \neq 0
\;,\qquad (C,s) \neq 0\; ;\label{B7.14a}\end{eqnarray}
the system (\ref{B7.6}) takes the form
$$
 0=0\;,\; 0=0\;,\; 0=0\;,\; 0= 0\;,
 \; 0=0\; ,\;  0=0\;, \; 0=0\;,\; 0=0\;,
 $$
 $$
 0 =0\;, \; 0=0\;, \; 0=0\;,\; 0=t+s\;, \;
  0=0\;, \;  0 =C+ D\;, \; s =C\;,\; t =D\;,
  $$
  $$
  0= 0\;,\; 0=0\;,\;  0=0\;,\; 0 =0\;,\;  0=0\;, \; 0 =0\;,\;  0=0\;, \;0=0\;,
  $$
  $$
  0=0\;,\; 0 =0\;,\ ;  0=0\;,\; 0=0\;,\ ;   C=s\;, \; D=t\;,\ ;
0=0\;, \; 0=0\;,
$$
$$
 0 =0\;, \; 0= 0\;,\;  0=t+C\;, \; 0=D +s\;,\;   0=0\;, \;0=0\;,\ ;
 0 =0\;, \; 0= D+C\;,
 $$
with the following solution
\begin{eqnarray}t=D\;,\; s=C =-D\;,\; (A,\alpha)=0\;,\;
    (B,\beta)=0\; ;\label{B7.14b}\end{eqnarray}
thus, we get

\vspace{3mm}

 \underline{solution $(KM-2)$},
\begin{eqnarray}{\bf n }=0\;,\; n_0=0\;,\; {\bf l }= D
({\bf m}-{\bf k})\;,\; l_0= D (m_0-k_0)\; ;\nonumber\\
G=\left(\ba{cc} K&0\\ D(M-K)&M
\ea\right) ;\qquad\qquad\nonumber\\
G'G=\left(\ba{cc} K'&0\\ D(M'-K')&M'
\ea\right) \left(\ba{cc} K'K&0\\
D(M'M-K'K)&M'M \ea\right).\label{B7.14d}\end{eqnarray}

Now, let it be (see (\ref{B7.9b}))
\begin{eqnarray}(A,\alpha)=0\;,\qquad (B,\beta) \neq 0\;,\qquad (D,t)=0
\;,\qquad (C,s) \neq 0\; ;\label{B7.15a}\end{eqnarray}
the system (\ref{B7.6}) reads
$$
 0=0\;,\; 0=0\;,\; \beta s=1\;,\; 0=0\;,\;
0= 0\;,\; \beta C=1\;,\; 0= 0\;,\; \beta=B\;,
$$
$$
 0 =0\;,\; 0 =0\;,\; 0=0\;,\; \beta s =1\;,\;
 0=0\;, \; s B =1\;,\; s=C\;,\;0 =0\;,
$$
$$
0=0\;,\;0=0\;,\;0=0\;,\; B =\beta\;,\;
 sB=1\;,\; BC=1\;,\ ;  0=0\;,\;  0 =0\;,
 $$
 $$
 0=0\;,\;B =B\;,\; B C=1\;,\;0=0\;,\;
  C=s\;, 0 =0\;,\ ; 0=0\;, \; 0=0\;,
 $$
 $$
  0= 0\;,\; \beta C= 1\;,\;  BC=1\;,\;  0=0\;,
\;0=0\;,\;  0 =0\;, \; BC=1\; ,
$$

\noindent which is equivalent to
\begin{eqnarray}s=C\;,\qquad\beta=B\;,\qquad BC=1\;.\label{B7.15b}\end{eqnarray}
Thus, we obtain

\vspace{3mm}

\underline{solution $(KM-3)$},
\begin{eqnarray}(A,\alpha)=0\;,\quad\beta=B\;,\quad (D,t)=0\;,
\quad C=s={1 \over B}\; ;\nonumber\\
 {\bf n }= B {\bf m}\;, \quad n_0=B m_0\;,\quad {\bf l }= {1 \over B}
{\bf k}\;,\;\;\quad l_0={1 \over B} k_0\; ;\label{B7.15c}
\end{eqnarray}
$$
G=\left(\ba{cc} K&B M\\ B^{-1} K&M \ea\right),\; G'G=
\left(\ba{cc} (K' +M')K&B (K' +M')M\\ B^{-1}(K'+M')K&(K'+M')M
\ea\right).
$$

Let us consider the case (see (\ref{B7.9b}))
\begin{eqnarray}(A,\alpha) \neq 0\;,\qquad (B,\beta)= 0\;,\qquad (D,t)=0\;,\qquad (C,s)
    \neq 0\; ;\label{B7.16a}\end{eqnarray}
the system (\ref{B7.6}) takes the form
$$
 \alpha s \alpha =0\;,\;0=0\;,\; 0=\alpha\;,\;0=0\;,\;
  0= 0\;,\;0=A\;,\; \alpha+\alpha AC=A\;,\;0=0\;,
  $$
  $$
  s \alpha s =0\;,\;0=0\;,\; 0 =0\;,\;0=s\;,
  \;0=0\;,\;0=C\;,\;s+ACs =C\;,\;0 =0\;,
  $$
  $$
  A \alpha C=0\;,\; 0=0\;,\; A+s A^2=\alpha\;,\;0  =0\;,
  \;  0=\alpha\;,\; 0 =A\;,\;  0= 0\;,\; 0=0\;,
  $$
  $$
A+AC A=A\;,\; 0 =0\;,\;  0= A\;,\;0=0\;,\;
 C+\alpha C^2=s\;,\; 0=0\;,\;  C s A=0\;,\; 0=0\;,
 $$
 $$
  0 =0\;,\; 0=0\;,\;  0=C\;,\; 0=s\;,\;
 AC^2 =0\;,\; 0=0\;,\;  0 =0\;,\;0=C\;,
 $$

 \noindent
which leads to already known

\vspace{3mm}

 \underline{solution $(KM-1)$},
\begin{eqnarray}(B,\beta)= 0\;,\; (D,t)=0\;,\; (A,\alpha)=0\;,\;(C,s)=0\;,\qquad
    G=\left(\ba {cc} K&0\\ 0&M\ea\right).\label{B.16c}\end{eqnarray}
Let us consider the case (see (\ref{B7.9b}))
\begin{eqnarray}(A,\alpha) \neq 0\;,\qquad (B,\beta) \neq 0\;,\qquad (D,t)= 0\;,\qquad
    (C,s)=0\;;\label{B7.17a}\end{eqnarray}
the system (\ref{B7.6}) gives
$$
0=0\;,\; 0=0\;,\;  0=\beta+\alpha\;,\; 0=0\;, \; 0=0\;,\;
0=B+A\;,\; \alpha=A\;,\; \beta=B\;,
$$
$$
 0 =0\;,\;0=0\;,\; 0=0\;,\;0=0\;,\;
 0= 0\;,\;0=0\;,\; 0=0\;,\;0 =D\;,
 $$
$$
0=0\;,\;0=0\;,\; A =\alpha\;,\;B=\beta\;,\; 0=B+\alpha\;,
 \;0 =A+\beta\;,\; 0=0\;,\;0=0\;,
 $$
$$
A=A\;,\;B =B\;,\; 0=A+ B\;,\;0=0\;,\; 0= 0\;,\;0 =0\;,\;
 0=0\;,\;0=0\;,
 $$
$$
 0 =0\;,\;0= 0\;,\; 0=0\;,\;0=0\;,\; 0=0\;,
0=0\;,\; 0 =0\;,\;0= 0\;,
 $$
which leads to

\vspace{2mm} \underline{solution $(KM-4)$},
\begin{eqnarray}
A= \alpha\;,\qquad B=\beta =-A\;,\qquad (D,t)=0\;,
\qquad (C,s)=0\;,\nonumber\\
 {\bf n }=A({\bf k}-{\bf m})\;,\; n_0=A(k_0-m_0)\;,\qquad {\bf l }= 0\;,\;\;l_0=0\;,\label{B7.17c}
\end{eqnarray}
$$
 G=\left(\ba{cc}K&A(K-M)\\0&M\ea\right),\;
 G'G=  \left(\ba{cc}K'K&A(K'K-M'M)\\0&M'M\ea\right).
 $$

Now, we are able to consider the four possibilities with three
vanishing blocks (see (\ref{B7.9c})). Let it be
\begin{eqnarray} (A,\alpha)=0\;,\qquad (B,\beta)\neq 0\;,\qquad (D,t) \neq 0
    \;,\qquad (C,s) \neq 0\;,\label{B7.18a}\end{eqnarray}
the system (\ref{B7.6}) gives
$$
 0=0\;,\; \beta t\beta=0\;,\;    \beta^2 s=\beta\;,\; 0=0\;,\;
    0=0\;,\; \beta C B= B\;,\;     0=0\;,\; \beta+\beta BD=B\;,
    $$
$$
0 =0\;,\; \beta t (s + t) =0\;,\;     st\beta =0\;,\; \beta s^2
=t+s\;,\;
     t CB=0\;,
     $$
     $$
     s BC =C+ D\;,\;     s =C\;,\; t+BD(t+s) =D\;,$$
$$
0=0\;,\; B\beta D =0\;,\;     0=0\;,\; B+tB^2 =\beta\;,\;
    sB^2=B\;,\; \beta BC =\beta\;,\;     0=0\;,\; 0=0\;,
    $$
$$
0= 0\;,\;  B+BD B =B\;,\;
    B^2 C=B\;,\; 0=0\;,\;
    C=s\;,\;D+\beta D^2+t BC =t\;,
    $$
    $$
        0=0\;,\; (t B +\beta C )=0\;, \;
 D s B =0\;,\; C\beta D=0\;,\;   \beta C^2=t+C\;,
 $$
 $$
 s BC=D +s\;,\;
    0=0\;,\; BD(D+C)=0\;,\;     DC B =0\;,\; BC^2=D+C\;.
    $$

These equations lead to the already known solution $(KM-1)$.

For the case
\begin{eqnarray}(A,\alpha) \neq 0\;,\qquad (B,\beta)= 0\;,\qquad (D,t) \neq
    0\;,\qquad (C,s) \neq 0\label{B7.19a}\end{eqnarray}
the system (\ref{B7.6}) gives only one and already known solution
$(KM-1)$.

Let it be
\begin{eqnarray}(A,\alpha) \neq 0\;,\qquad (B,\beta)\neq 0\;,\qquad (D,t)=0\;,
    \qquad (C,s) \neq 0\; ;\label{B7.20a}\end{eqnarray}
the system (\ref{B7.6})  leads to the already known solution
$(KM-1)$.

Consider the  case where we assume
\begin{eqnarray}(A,\alpha) \neq 0\;,\qquad (B,\beta)\neq 0\;,\qquad (D,t)
\neq 0\;,\qquad (C,s)=0\; ;\label{B7.21a}\end{eqnarray}
the system (\ref{B7.6}) reads
$$
 0= 0\;,\;  \beta t (\alpha+\beta)=0\;,\;
    \alpha^2 t=\beta+\alpha\;,\; \alpha\beta t=0\;,
    $$
    $$
    \beta AD=0\;,\; \alpha AD=B+A\;,\;     \alpha=A\;,\; \beta+(\alpha+\beta) BD=B\;,
    $$
    $$
     0=0\;,\;  \beta t t =0\;,\;     0=0\;,\; \alpha t^2=t\;,
     $$
     $$
    0=0\;,\; t AD=D\;,\;
    0=0\;,\; t+BD t =D\;,$$
$$
0=0\;, \; B (\beta D +t A)=0\;,\;     A =\alpha\;,\; B+tB^2 +\beta
AD=\beta\;,
$$
$$
    \alpha AD=B+\alpha\;,\;  \beta tA^2 =A+\beta\;,\;
    t AB=0\;,\; \alpha BD =0\;, $$
$$
A=A\;,\; B+BD (B+A) =B\;,\;     A^2 D=A+ B\;,\; AB D=0\;,
$$
$$
    0= 0\;,\; D+\beta D^2 =t\;,\;
    0=0\;,\;D t B =0\;, $$
$$
 0 =0\;,\; 0=0\;,\;     t AD=t\;,\; \alpha D^2 =D\;,
 $$
 $$
    0=0\;,\; BD D =0\;,\;
    0 =0\;,\; AD^2=D\;.
    $$
    which leads to the already known solution $(KM-1)$.

It remains to examine the case of four non-zero blocks
\begin{eqnarray}(A,\alpha) \neq 0\;,\qquad (B,\beta) \neq 0\;,\qquad (D,t) \neq 0
    \;,\qquad (C,s) \neq 0\;.\label{B7.22a}\end{eqnarray}
The system (\ref{B7.6}) can be simplified to
$$\left.\ba{ll} (\alpha+\beta)=0\;,\qquad&(\alpha+\beta)=0\;,\\[2mm]
    \alpha^2 t+\beta^2 s=\beta+\alpha\;,\qquad&(s+ t)=0\;,\\[2mm]
    \alpha BC+\beta AD=0\;,\qquad&\alpha AD+\beta BC=B+A\;,\\[2mm]
    \alpha+(\alpha +\beta) AC=A\;,\qquad&\beta+(\alpha
    +\beta) BD=B\;,\ea \right.$$
$$\left.\ba{ll} (s +t) =0\;,\qquad&(s + t) =0\;,\\[2mm]
    (\alpha+\beta)=0\;,\qquad&\beta s^2+\alpha t^2=t+s\;,\\[2mm]
    sAD+t CB=0,\qquad&s BC+t AD=C+ D\;,\\[2mm]
    s+AC(s+t)=C\;,\qquad&t+BD(t+s) =D\;,\ea \right.$$
$$ \left.\ba{ll}(\alpha C+s B)=0\;,\qquad&(\beta D +t A)=0\;,\\[2mm]
    A+s A^2+\alpha BC=\alpha\;,\qquad&B+tB^2 +\beta AD=\beta\;,\\[2mm]
    \alpha AD+sB^2=B+\alpha\;,\qquad&\beta BC+tA^2 =A+\beta\;,\\[2mm]
  \beta C+t B=0\;,\qquad&\alpha D+s A =0\;, \ea \right.$$
$$\left.\ba{ll} (A+B)=0\;,\qquad&(B+A) =0\;,\\[2mm]
    A^2 D+B^2 C=A+ B\;,\qquad&(C+D)=0\;,\\[2mm]
    C+\alpha C^2+s AD=s\;,\qquad&D+\beta D^2+t BC =t\;,\\[2mm]
    (s A+\alpha D)=0\;,\qquad&(t B +\beta C )=0\;, \ea \right.$$
\begin{eqnarray}\left.\ba{ll} ( \alpha C+s B) =0\;,\qquad&(\beta D+t A)=0\;,\\[2mm]
  \beta C^2+t AD=t+C\;,\qquad&\alpha D^2 +s BC=D +s\;,\\[2mm]
    (C+D)=0\;,\qquad&(D+C)=0\;,\\[2mm]
    (A+B) =0\;,\qquad&BC^2+AD^2=D+C\;. \ea \right.\label{B7.22b}\end{eqnarray}
By eliminating the variables
\begin{eqnarray}B =-A\;,\;\beta =-\alpha\;,\qquad C =-D\;,\; s=-t\;,\label{B7.22c}\end{eqnarray}
we get the system
$$
 0=0\;,\; 0=0\;,\;    0= 0\;,\; 0=0\;,
 $$
 $$
    0=0\;,\; 0= 0\;,\;     \alpha=A\;,\; \beta=-A\;,
    $$
$$
0 =0\;,\; 0=0\;,\;     0=0\;,\;0=0\;,
$$
$$
    0=0\;,\;  0=0\;,\;
    s =-D\;,\; t =D\;,
    $$
$$ \left.\ba{ll}(-\alpha D-s A)=0\;,\qquad&(\beta D +t A)=0\;,\\[2mm]
    A+s A^2+\alpha AD=\alpha\;,\qquad&-A+tA^2 +\beta AD=\beta\;,\\[2mm]
    \alpha AD+sA^2=-A+\alpha\;,\qquad&\beta AD+tA^2 =A+\beta\;,\\[2mm]
  \beta -D-t A=0\;,\qquad&\alpha D+s A =0\;, \ea \right.$$
$$\left.\ba{ll} 0=0\;,\qquad&0 =0\;,\\[2mm]
    0=0\;,\qquad&0=0\;,\\[2mm]
    -D+\alpha D^2+s AD=s\;,\qquad&D+\beta D^2+t AD =t\;,\\[2mm]
    (s A+\alpha D)=0\;,\qquad&(-t A -\beta D )=0\;, \ea \right.$$
$$\left.\ba{ll} ( \alpha C-s A) =0\;,\qquad&(\beta D+t A)=0\;,\\[2mm]
  \beta C^2+t AD=t+C\;,\qquad&\alpha D^2 -s AC=D +s\;,\\[2mm]
    0=0\;,\qquad&0=0\;,\\[2mm]
    0 =0\;,\qquad&0= 0\;. \ea \right.$$
From this it follows
\begin{eqnarray}A=\alpha\;,\qquad B=\beta=-A\;,\qquad t=D\;,\qquad
    C=s=-D\;,\label{B7.22d}\end{eqnarray}
which describes the following

\vspace{3mm}

\underline{solution $(KM-5)$},
\begin{eqnarray} {\bf n }=A ({\bf k}-{\bf m})\;,\qquad n_0=A (k_0-m_0)\;,\nonumber\\
{\bf l }= D ({\bf m} -{\bf k})\;,\;\;\qquad l_0=D (m_0-k_0)\; ;\nonumber\label{B7.23a}\\
 G=\left(\ba{cc}K&A(K-M)\\C(K-M)&M\ea\right).\label{B7.23b} \end{eqnarray}

Let us verify by direct calculation that the matrices with such
structure form a group:
\begin{eqnarray}\left(\ba{cc}
K'&A (K'-M')\\
C(K'-M')&M' \ea\right) \left(\ba{cc}
K&A (K-M)\\
C(K-M)&M
\ea\right)\nonumber\\
=\left(\ba{cc}
K''&A (K''-M'')\\
C(K''-M'')&M'' \ea\right).\qquad\qquad \nonumber\end{eqnarray}
\section{Two Independent Vectors: Variant II(l,n)}
Let us examine {\bf variant II(l,n)}:
\begin{eqnarray} {\bf k}=(A {\bf l}+B {\bf n})\;,
\qquad k_0=(\alpha l_0+\beta n_0)\;,\nonumber\\
{\bf m} =(D {\bf l}+C {\bf n})\;,\qquad
 m_0=(t l_0+s n_0)\;.\label{B8.1}\end{eqnarray}
The multiplication law takes the form
$$
k_0''=(\alpha l_0'+\beta n_0')\; (\alpha l_0+\beta n_0)+(A {\bf
l}'+B {\bf n}') (A {\bf l}+B {\bf n}) +n'_0\; l_0+{\bf n}'\cdot
{\bf l}\;,
$$
$$
m_0''=(t l_0'+s n_0')\; (t l_0+s n_0) + (D {\bf l}'+C {\bf n} ')\;
(D {\bf l}+C {\bf n}) +l'_0\; n_0+{\bf l}'\cdot {\bf n}\;,
$$
$$
n_0''=(\alpha l_0'+\beta n_0')\; n_0+(A {\bf l}'+B {\bf n}')\;
{\bf n} +n'_0\; (t l_0+s n_0)+{\bf n}'\; (D {\bf l}+C {\bf n})\;,
$$
$$
l_0''=l_0'\; (\alpha l_0+\beta n_0)+{\bf l}'\; (A {\bf l}+B {\bf
n}) +(t l_0 '+s n_0')\; l_0+(D {\bf l} '+C {\bf n} ')\; {\bf l}\;,
 $$

$$ {\bf k}''=(\alpha l_0'+\beta n_0')\; (A {\bf l}+B
{\bf n})+(A {\bf l}'+B {\bf n}')\; (\alpha l_0+\beta n_0)$$
$$+ i\; (A {\bf l}'+B {\bf n}')\times (A {\bf l}+B {\bf n})
+n_0'\; {\bf l}+{\bf n}'\; l_0+i\; {\bf n}'\times {\bf l}\;,
$$
$$
{\bf m}''=(t l_0'+s n_0')\; (D {\bf l}+C {\bf n}) +
 (D {\bf l}'+C {\bf n}')\; (t l_0+s n_0)
 $$
 $$
+i\; (D {\bf l} '+C {\bf n} ')\times (D {\bf l}+C {\bf n})+l_0'\;
{\bf n}+{\bf l}'\; n_0+i\; {\bf l}' \times {\bf n}\;,
$$
$$
{\bf n}''=(\alpha l_0'+\beta n_0')\; {\bf n}+(A {\bf l}'+B {\bf
n}')\; n_0 +
 i\; (A {\bf l} '+ B {\bf n}')\times {\bf n}
$$
$$
+n_0'\; (D {\bf l}+C {\bf n}) +{\bf n}'\; (t l_0+s n_0)+i\; {\bf
n}'\times (D {\bf l}+C {\bf n})\;,
$$
\begin{eqnarray}{\bf l}''=l'_0\; (A {\bf l}+B {\bf n})+{\bf l}'\;
(\alpha l_0+\beta n_0)+i\; {\bf l}'\times (A {\bf l}
+ B {\bf n})\nonumber\\
+ (t l_0'+s n_0')\; {\bf l}+(D {\bf l}'+C {\bf n} ')\; l_0+i\;
 (D {\bf l} '+C {\bf n}')\times {\bf l }\;.\label{B8.2}\end{eqnarray}

The equation
\begin{eqnarray}{\bf k}''=A {\bf l}''+B {\bf n} ''\label{B8.3a}\end{eqnarray}
takes the form
\begin{eqnarray} (\alpha l_0'+\beta n_0') (A {\bf l}+B {\bf n}) +
(A {\bf l}'+B {\bf n}') (\alpha l_0+\beta n_0)\nonumber\\
+
 i (A {\bf l}'+B {\bf n}')\times (A {\bf l}+B {\bf n})
+n_0' {\bf l}+{\bf n}' l_0+i
{\bf n}'\times {\bf l}\nonumber\\
= A\; [ l'_0 (A {\bf l}+B {\bf n})+{\bf l}' (\alpha l_0
+\beta n_0)+i {\bf l}'\times (A {\bf l}+B {\bf n})\nonumber\\
 +
(t l_0'+s n_0') {\bf l}+(D {\bf l}'+C {\bf n} ') l_0+i
 (D {\bf l} '+C {\bf n}')\times {\bf l }]
\nonumber
\\
 +
B\; [ (\alpha l_0'+\beta n_0') {\bf n}+(A {\bf l}' + B {\bf n}')
n_0 +
 i (A {\bf l} '+ B {\bf n}')\times {\bf n}
\nonumber
\\
+n_0' (D {\bf l}+C {\bf n}) +{\bf n}'\; (t l_0+s n_0)+i {\bf
n}'\times (D {\bf l}+C {\bf n}) ]\; ; \nonumber\end{eqnarray}
whence it follows
\begin{eqnarray}l'_0 {\bf ll}:\qquad \alpha A =A^2+At\;,\quad
n'_0 {\bf n}:\qquad\beta B=\beta B+BC\;,\nonumber\\
n'_0 {\bf l}:\qquad\beta A +1=As +BD\;,\quad
l'_0 {\bf n}:\qquad \alpha B =AB+\alpha B\;,\nonumber\\
n_0 {\bf l}' :\qquad\beta A=\beta A +BA\;,\quad
l_0 {\bf n}' :\qquad \alpha B +1=AC+B t\;,\nonumber\\
l_0 {\bf l}' :\qquad \alpha A=\alpha A+AD\;,\quad n_0 {\bf n}'
:\qquad\beta B=B^2+B s\;,
\nonumber\\
{\bf l}'\times {\bf l} :\qquad A^2=A^2 +AD\;,\quad
{\bf n}'\times {\bf n} :\qquad B^2=B^2+BC\;,\nonumber\\
{\bf n}'\times {\bf l} :\qquad BA+1=AC+BD\;,\quad {\bf l}'\times
{\bf n} :\qquad AB=AB +AB\,.\label{B8.3b}\end{eqnarray}

The equation
\begin{eqnarray}{\bf m}''=D {\bf l}'' +C {\bf n} ''\label{B8.4a}\end{eqnarray}
takes the form
\begin{eqnarray} (t l_0'+s n_0') (D {\bf l}+C {\bf n}) +
 (D {\bf l}'+C {\bf n}') (t l_0+s n_0)\nonumber\\
+ i (D {\bf l} '+C {\bf n} ')\times (D {\bf l}+C
{\bf n})+l_0' {\bf n}+{\bf l}' n_0+i\; {\bf l}'\cdot {\bf n} =\nonumber\\
\times D\; [ l'_0 (A {\bf l}+B {\bf n})+{\bf l}' (\alpha
l_0+\beta n_0)+i {\bf l}'\times (A {\bf l}+B {\bf n})\nonumber\\
+ (t l_0'+s n_0') {\bf l}+(D {\bf l}'+C {\bf n} ') l_0+i
 (D {\bf l} '+C {\bf n}')\times {\bf l }] \nonumber \\
 + C\; [ \alpha l_0'+\beta n_0')\; {\bf n}+(A {\bf l}' + B {\bf n}')\; n_0 +
 i\; (A {\bf l} '+ B {\bf n}')\times {\bf n} \nonumber \\
+n_0'\; (D {\bf l}+C {\bf n}) +{\bf n}'\; (t l_0+s n_0)+i\; {\bf
n}'\times (D {\bf l}+C {\bf n})\;,\nonumber\end{eqnarray}
from which we obtain
$$ l'_0 {\bf l}:\qquad tD=DA+t D\;,\quad
 n'_0 {\bf n}:\qquad s C=C\beta +C^2\;,
 $$
$$
n'_0 {\bf l}:\qquad s D=sD+CD\;,\quad l'_0 {\bf n}:\qquad
tC+1=DB+\alpha C\;,$$
$$
n_0 {\bf l}' :\qquad sD +1=\beta D+CA\;,\quad l_0 {\bf n}' :\qquad
t C=DC+t C\;,
$$
$$
l_0 {\bf l}' :\qquad t D=\alpha D +D^2\;,\quad n_0 {\bf n}'
:\qquad sC=CB+sC\;,
$$
$$
{\bf l}'\times {\bf l} :\qquad D^2=DA+D^2\;,\quad {\bf n}'\times
{\bf n} :\qquad C^2=CB+C^2\;,
$$
$$
{\bf n}'\times {\bf l} :\qquad CD=CD+CD\;,\quad {\bf l}'\times
{\bf n} :\qquad CD+1=DB+CA\;.
$$

The equation
\begin{eqnarray}k_0''=\alpha l_0''+\beta n_0 ''\label{B8.5a}\end{eqnarray}
gives
\begin{eqnarray} (\alpha l_0'+\beta n_0')\; (\alpha l_0+\beta n_0) +
(A {\bf l}'+B {\bf n}') (A {\bf l}+B {\bf n}) +n'_0\;
l_0+{\bf n}'\cdot {\bf l}\nonumber\\
= \alpha\; [\; l_0'\; (\alpha l_0+\beta n_0)+{\bf l}'\; (A {\bf
l}+B {\bf n})
+(t l_0 '+s n_0')\; l_0+(D {\bf l} '+C {\bf n} ')\;{\bf l} ]\nonumber\\
+\beta\; [\; (\alpha l_0'+\beta n_0')\; n_0+(A {\bf l}'+B {\bf
n}')\; {\bf n} +n'_0\; (t l_0+s n_0)+{\bf n}'\; (D {\bf l}+C {\bf
n}  ) ]\;,\nonumber \end{eqnarray}
which infers
$$
l_0'l_0:\qquad \alpha^2=\alpha^2+\alpha t\;,\quad
n_0'n_0:\qquad\beta^2=\beta^2+s\beta\;,
$$
$$
l_0'n_0:\qquad \alpha\beta=\alpha\beta+\alpha\beta\; ,\quad
n_0'l_0:\qquad \alpha\beta +1=\alpha s+\beta t\;,
$$
$$
{\bf l}' \cdot{\bf l} :\qquad A^2=\alpha A+\alpha D\;,\quad {\bf
n}'\cdot {\bf n} :\qquad B^2=\beta B+\beta C\;,
$$
$$
{\bf l}'\cdot {\bf n} :\qquad AB=\alpha B+\beta A\;,\quad {\bf
n}'\cdot {\bf l} :\qquad AB +1=\alpha C+\beta D\;.
$$

The equation
\begin{eqnarray}m_0 ''=t l_0 ''+s n_0 ''\label{B8.6a}\end{eqnarray}
takes the form
\begin{eqnarray} (t l_0'+s n_0')\; (t l_0+s n_0) +
(D {\bf l}'+C {\bf n} ')\; (D {\bf l}+C {\bf n})
+l'_0\; n_0+{\bf l}'\cdot {\bf n}\nonumber\\
=t\; [\; l_0'\; (\alpha l_0+\beta n_0)+{\bf l}'\; (A {\bf l}+B
{\bf n})
+(t l_0 '+s n_0')\; l_0+(D {\bf l} '+C {\bf n} ')\; {\bf l }\; ]\nonumber\\
 + s\; [\;(\alpha l_0'+\beta n_0')\; n_0+(A {\bf l}'+B {\bf n}')\; {\bf n}
+n'_0\; (t l_0+s n_0)+{\bf n}'\; (D {\bf l}+C {\bf
n})\;]\;,\nonumber \end{eqnarray}
from which we obtain
$$l_0'l_0:\qquad t^2=t \alpha+t^2\;,\quad
 n_0'n_0:\qquad s^2=s\beta+s^2\;,
 $$
$$
l_0'n_0:\qquad ts +1=t\beta+s \alpha\;,\quad
 n_0'l_0:\qquad st=st+st\;,
 $$
$${\bf l}' \cdot{\bf l} :\qquad D^2=tA+tD\;,\quad
{\bf n}'\cdot {\bf n} :\qquad C^2= sB+s C\;,$$
$$
{\bf l}'\cdot {\bf n} :\qquad DC +1=tB+sA\;,\quad
    {\bf n}'\cdot {\bf l} :\qquad CD= tC+ sD\;.
    $$

We further collect the equations together, and infer
\begin{eqnarray} \alpha A =A^2+At\;,\; 0= BC\;,\;
\beta A +1=As +BD\;,\;  0 =AB\;,\;  0=+BA\;,
\nonumber\end{eqnarray}
\begin{eqnarray} \alpha B +1=AC+B t\;,\;  0= AD\;,\;
\beta B=B^2+B s\;,\;  0=AD\;,\;  0=BC\;, \nonumber\end{eqnarray}
\begin{eqnarray}
 BA+1=AC+BD\;,\;
 0=AB\;,\;
 0=DA\;,\;
 s C=C\beta +C^2\;,\;
 0= CD\;,\nonumber
 \end{eqnarray}
 \begin{eqnarray}
  tC+1=DB+\alpha C\;,\;
  sD +1=\beta D+CA\;,\;
  0= DC\;,\;
  t D=\alpha D +D^2\;,\;
   0= CB\;,\nonumber
   \end{eqnarray}
   \begin{eqnarray}
   0=DA\;,\;
    0=CB\;,\;
 0=CD\;,\;
 CD+1=DB+CA\; ;
\nonumber\end{eqnarray}
\begin{eqnarray}0=\alpha t\;,\quad
 0=s\beta\;,\quad
 0=\alpha\beta\;,\quad
 \alpha\beta +1=\alpha s+\beta t\;,
 \nonumber
 \end{eqnarray}
 \begin{eqnarray}
  A^2=\alpha A+\alpha D\;,\;
 B^2=\beta B+\beta C\;,\;
 AB=\alpha B+\beta A\;,\;
  AB +1=\alpha C+\beta D\; ;
\nonumber\end{eqnarray}
\begin{eqnarray} 0=t \alpha\;,\quad
0=s\beta\;,\quad
 ts +1=t\beta+s \alpha\;,\quad
 0=st\;,
\nonumber\end{eqnarray}
\begin{eqnarray} D^2=tA+tD\;,\;
 C^2= sB+s C\;,\;
 DC +1=tB+sA\;,\;
 CD= tC+ sD\;.
\nonumber\end{eqnarray}

From the very beginning we  note that the trivial set of all
vanishing parameters  provide no
 solution of the system.

Let us separate the most simple equations
\begin{eqnarray}BC=0\;,\quad\beta s =0\;,\quad
AB=0\;,\quad \alpha\beta =0\;,\nonumber\\
 AD=0\;,\quad \alpha t =0\;,\quad
 CD =0\;,\quad st =0\; ;\label{B8.7a}\end{eqnarray}
they can be satisfied by six solutions
\begin{eqnarray}1)\quad B,\beta=0\;,\quad D,t =0\; ;\quad
2)\quad A,\alpha=0\;,\quad C, s=0\; ;\nonumber\\
3)\qquad \,B,\beta=0\;,\qquad A, \alpha =0\;,\qquad C, s=0\; ;\nonumber\\
4)\qquad \, A,\alpha=0\;,\qquad B,\beta =0\;,\qquad D,t=0
\; ;\nonumber\\
5)\qquad C,s=0\;,\qquad \; B,\beta =0\;,\qquad D,t=0\; ; \nonumber\\
6)\qquad D,t =0\;,\qquad\;\;  C,s =0\;,\qquad A, \alpha =0\;.
\label{B8.7b}\end{eqnarray}

The remaining equations are
\begin{eqnarray} \alpha A =A^2+At\;,\quad
\beta A +1=As +BD\;,\quad
 \alpha B +1=AC+B t\;,\nonumber\\
\beta B=B^2+B s\;,\quad
 BA+1=AC+BD\; ;\label{B8.8a}\end{eqnarray}
\begin{eqnarray} s C=C\beta +C^2\;,\quad
 tC+1=DB+\alpha C\;,\quad
 sD +1=\beta D+CA\;,\nonumber\\
 t D=\alpha D +D^2\;,\quad
 CD+1=DB+CA\,;\label{B8.8b}\end{eqnarray}
\begin{eqnarray} \alpha\beta +1=\alpha s+\beta t\;,\quad
 A^2=\alpha A+\alpha D\;,\quad
 B^2=\beta B+\beta C\;,\nonumber\\
 AB=\alpha B+\beta A\;,\quad
 AB +1=\alpha C+\beta D\; ;\label{B8.8c}\end{eqnarray}
\begin{eqnarray} ts +1=t\beta+s \alpha\;,\quad
 D^2=tA+tD\;,\quad
 C^2= sB+s C\;,\nonumber\\
 DC +1=tB+sA\;,\quad
 CD= tC+ sD\;.\label{B8.8d}\end{eqnarray}

We see that variants 3), 4), 5), 6) are not consistent with the
additional equations
    (\ref{B8.8a})--(\ref{B8.8d}). Therefore, we shall further examine only the variants 1) and 2).

In the case 1), we have
\begin{eqnarray}1)\qquad B,\beta=0\;,\qquad D,t =0\; ;\nonumber\end{eqnarray}
the equations (\ref{B8.8a})--(\ref{B8.8d}) give
\begin{eqnarray} \alpha A =A^2\;,\qquad  1=As\;,\qquad 1=AC\;,\nonumber\\
 s C=C^2\;,\qquad 1=\alpha C\;, \qquad 1=CA\;,\nonumber\\
 1=\alpha s\;,\qquad  A^2=\alpha A\;,\qquad 1=\alpha C\;,\nonumber\\
 1= s \alpha\;,\qquad  C^2= s C\;,\qquad  1=sA\; ; \nonumber\end{eqnarray}
the last system has only one solution
\begin{eqnarray}B,\beta=0\;,\qquad D,t =0\;,\qquad \alpha=A \neq 0\;,
    \qquad C=s={1 \over A}\; ;\label{B8.10a}\end{eqnarray}
so we arrive at

\vspace{3mm}

\underline{solution $(LN-1)$},
\begin{eqnarray} {\bf k}=A {\bf l}\;,\; k_0=A l_0\;,\; {\bf m}=A^{-1} {\bf n}\;,\;
 m_0=A^{-1} n_0\;,\nonumber\label{B8.10b}\\
G=\left(\ba{cc}A L&N\\L&A^{-1} N\ea\right),\qquad\qquad\qquad\label{B8.10c}\\
G'G
= \left(\ba{cc}A (AL'L +A^{-1} N'L)&(AL'N+A^{-1} N' N)\\
(AL'L +A^{-1} N'L)&A^{-1} (AL'N+A^{-1} N'
N)\ea\right).\nonumber\end{eqnarray}
This is a set of degenerate matrices of rank 2.

For the case
\begin{eqnarray}2)\qquad A,\alpha=0\;,\qquad C, s=0\;
\nonumber\end{eqnarray}
the equations (\ref{B8.8a})--(\ref{B8.8d}) give
\begin{eqnarray} 1=BD\;,\qquad 1=B t\;,\qquad\beta B=B^2\;,\nonumber\\
 1=DB\;,\qquad 1=\beta D\;,\qquad t D=D^2\;,\nonumber\\
 1=\beta t\;,\qquad B^2=\beta B\;,\qquad 1=\beta D\;,\nonumber\\
 1=t\beta\;,\qquad D^2=tD\;, \qquad 1=tB\;,\nonumber\end{eqnarray}
which has only one solution
\begin{eqnarray}A,\alpha=0\;,\qquad C, s=0\;,\qquad\beta=\beta \neq
0\;,\qquad D=t=B^{-1}\; ;\label{B8.12b}\end{eqnarray}
so we obtain

\vspace{3mm}

\underline{solution $(LN-2)$},
\begin{eqnarray} {\bf k}=B {\bf n}\;,\; k_0=B n_0\;,\; {\bf m} =B^{-1} {\bf l}\;,\;
 m_0=B^{-1} l_0\;,\nonumber\\
G=\left(\ba{cc}B N&N\\L&B^{-1} L\ea\right),\qquad\qquad\qquad\label{B8.12c}\\
G'G=  \left(\ba{cc} B
(BN'N+B^{-1}N'L)&(BN'N+B^{-1}N'L)\\(BL'N+B^{-1}L'L)&B^{-1}
(BL'N+B^{-1}L'L) \ea\right).\nonumber\end{eqnarray}
This is a set of degenerate matrices of rank 2.
The analysis of the {\bf variant} {\bf II(l,n)} is completed.
\section{Two Independent Vectors: Variant II(k,l)}
Let us examine the {\bf variant II(k,l)}:
\begin{eqnarray} {\bf n}=A {\bf k}+B {\bf l}\;,\qquad n_0=\alpha k_0+\beta l_0\;,\nonumber\\
{\bf m}=D {\bf l}+C {\bf k}\;,\qquad m_0=t l_0+s
k_0\;.\label{B9.1}\end{eqnarray}
The multiplication law takes the form
\begin{eqnarray}k_0''=k_0'\; k_0+{\bf k}'\cdot {\bf k}+(\alpha k_0' +
\beta l_0')\; l_0+(A {\bf k}'+B {\bf l}')\; {\bf l}\;,\nonumber\\
m_0''=(t l_0'+s k_0') (t l_0+s k_0)+(D{\bf l}'+C {\bf k}')\; (D {\bf l}+C {\bf k})\nonumber\\
+ l'_0\; (\alpha k_0+\beta l_0)+{\bf l}' (A {\bf k}+ B {\bf l})\;,\nonumber\\
n_0''=k_0'\; (\alpha k_0+\beta l_0)+{\bf k}'\;(A {\bf k}+B {\bf l})\nonumber\\
 + (\alpha k_0'+\beta l_0')\; (t l_0+s k_0)+(A {\bf k}'+B {\bf l}')\;
 (D {\bf l}+C {\bf k})\;,\nonumber\\
l_0''=l_0'\; k_0+{\bf l}'\cdot {\bf k} + (t l_0'+s k_0') l_0+(D
{\bf l}'+ C {\bf k}')\; {\bf l}\;,\nonumber\end{eqnarray}
\begin{eqnarray}{\bf k}''=k'_0\; {\bf k}+{\bf k}'\; k_0+i\; {\bf
k}'\times {\bf k} +(\alpha k_0'+\beta l_0')\; {\bf l}+(A {\bf k}'\nonumber\\
+B{\bf l}')\; l_0+i\; (A {\bf k}'+B {\bf l}')\times {\bf l}\;,\nonumber\\
{\bf m}''=(t l_0'+s k_0') (D {\bf l}+C {\bf k}) +(D {\bf l} '+C {\bf k}')\; (t l_0+s k_0)\nonumber\\
+i\;(D {\bf l}'+C {\bf k}')\times (D {\bf l}+C {\bf k})\nonumber\\
 + l_0'\; (A {\bf k}+B {\bf l})+{\bf l}'\; (\alpha k_0+\beta l_0)+i\;
{\bf l}'\times (A {\bf k}+B {\bf l})\;,\nonumber\\
{\bf n}''=k'_0\; (A {\bf k}+B {\bf l})+{\bf k}'\;
(\alpha k_0+\beta l_0) +i\; {\bf k}'\times (A {\bf k}+B {\bf l})\nonumber\\
 + (\alpha k_0'+\beta l_0')\; (D {\bf l}+C {\bf k})\nonumber\\
+(A {\bf k}'+B {\bf l}')\; (t l_0+s k_0)+i\;
(A {\bf k}'+B {\bf l}')\times (D {\bf l}+C {\bf k})\;,\nonumber\\
{\bf l}''=l'_0\; {\bf k}+{\bf l}'\; k_0+i\; {\bf l}'\times {\bf k}\nonumber\\
 + (t l_0'+s k_0') {\bf l} + (D {\bf l}'+C {\bf k}')\; l_0+i\;
(D {\bf l}'+C {\bf k}')\times {\bf l
}\;.\label{B9.2}\end{eqnarray}

The equation
\begin{eqnarray} {\bf n}''=A {\bf k}''+B {\bf l}''\label{B9.4a}\end{eqnarray}
gives
\begin{eqnarray} k'_0\; (A {\bf k}+B {\bf l}) +
{\bf k}'\; (\alpha k_0+\beta l_0) +i\; {\bf k}'\times (A {\bf k}+B {\bf l})\nonumber\\
 + (\alpha k_0'+\beta l_0')\; (D {\bf l}+C {\bf k})+(A {\bf k}'+B {\bf l}')\; (t l_0+s k_0)+i\;
(A {\bf k}'+B {\bf l}')\times (D {\bf l}+C {\bf k})\nonumber\\
= A\; [\; k'_0\; {\bf k}+{\bf k}'\; k_0+i\; {\bf k}'\times {\bf
k}+(\alpha k_0'+\beta l_0')\; {\bf l} +
(A {\bf k}'+B {\bf l}')\; l_0+i\; (A {\bf k}'+B {\bf l}')\times {\bf l}\; ]\nonumber\\
+ B\; [\; l'_0\; {\bf k}+{\bf l}'\; k_0+i\; {\bf l}'\times {\bf k}
+
 (t l_0'+s k_0') {\bf l} + (D {\bf l}'+C {\bf k}')\; l_0+i\;
(D {\bf l}'+C {\bf k}')\times {\bf l
}\;].\nonumber\\\label{B9.4b}\end{eqnarray}

Further, we obtain the system
$$
k'_0\; {\bf k}:\qquad \alpha C=0\;,\quad k_0\; {\bf k}':\qquad
\alpha+A s=A\;,
$$
$$
k'_0\; {\bf l}:\qquad B+\alpha D=A \alpha+B s\;,\quad k_0\; {\bf
l}':\qquad sB=B\;,
$$
$$
l'_0\; {\bf l}:\qquad\beta D=A\beta+B t\;,\quad l_0\; {\bf
l}':\qquad B t=A B+BD\;,
$$
$$
l'_0\; {\bf k}:\qquad\beta C=B\;,\quad l_0\; {\bf k}':\qquad\beta
+At=A^2+BC\;,
$$
$$
{\bf k}'\times {\bf k} :\qquad AC=0\;,\quad {\bf l}'\times {\bf l}
:\qquad AB=0\;,\
$$
$$
{\bf k} '\times {\bf l} :\qquad B+AD=A^2+BC\;,\quad {\bf l'}\times
{\bf k} :\qquad BC=B\;.
$$

The equation
\begin{eqnarray}n_0''=\alpha k_0''+\beta l_0''
\nonumber \label{B9.5a}\end{eqnarray} leads to
\begin{eqnarray}k_0'\; (\alpha k_0+\beta l_0)+{\bf k}'\; (A {\bf k}+ B {\bf l})\nonumber\\
 + (\alpha k_0'+\beta l_0')\; (t l_0+s k_0)+(A {\bf k}'+B {\bf l}')\; (D {\bf l}+C {\bf k}) =\nonumber\\
\alpha\; [\; k_0'\; k_0+{\bf k}'\; {\bf k}+(\alpha k_0'+\beta
l_0')\; l_0+(A {\bf k}'+B {\bf l}')
\; {\bf l}\; ]\nonumber\\
 +\beta\; [\; l_0'\; k_0+{\bf l}'\cdot {\bf k} +
 (t l_0'+s k_0') l_0+(D {\bf l}'+C {\bf k}')\cdot {\bf l}\; ]\; ;\nonumber \end{eqnarray}
this infers
$$
k'_0k_0:\qquad \alpha s=0\;,\quad l'_0l_0:\qquad 0=\alpha\beta\;,
$$
$$
k'_0l_0:\qquad\beta+\alpha t=\alpha^2+\beta s\;,\quad
l'_0k_0:\qquad\beta s=\beta\;,
$$
$$
{\bf k}'\cdot {\bf k}:\qquad A+AC=\alpha\;,\quad {\bf l}'
\cdot{\bf l}:\qquad BD=\alpha B+\beta D\;,
$$
$$
{\bf k}' \cdot{\bf l}:\qquad B +AD=\alpha A+\beta C\;,\quad {\bf
k} \cdot{\bf l}':\qquad BC=\beta\;.
$$

The equation
\begin{eqnarray}{\bf m}''=D {\bf l}''+C {\bf k}''\label{B9.6a}\end{eqnarray}
leads to
\begin{eqnarray} (t l_0'+s k_0') (D {\bf l}+C {\bf k})+(D {\bf l} '+C {\bf k}')\;
(t l_0+s k_0)+i\; (D {\bf l}'+C {\bf k}')
\times (D {\bf l}+C {\bf k})\nonumber\\
 + l_0'\; (A {\bf k}+B {\bf l})+{\bf l}'\; (\alpha k_0+\beta l_0)+i\;
{\bf l}'\times (A {\bf k}+B {\bf l}) = D\; [\; l'_0\; {\bf k}+{\bf
l}'\; k_0+i\; {\bf l}'
\times {\bf k}\nonumber\\
+ (t l_0'+s k_0') {\bf l} + (D {\bf l}'+C {\bf k}')\; l_0+i\; (D
{\bf l}'+C {\bf k}')\times {\bf l }\;]
+ C\; [\; k'_0\; {\bf k}+{\bf k}'\; k_0+i\; {\bf k}'\times {\bf k}\nonumber\\
 +(\alpha k_0'+\beta l_0')\; {\bf l}+(A {\bf k}'+B
{\bf l}')\; l_0+i\; (A {\bf k}'+B {\bf l}')\times {\bf l}\; ]\;
;\nonumber\end{eqnarray}
further we obtain
$$
k'_0\; {\bf k}:\qquad sC=C\;,\quad k_0\; {\bf k}':\qquad sC=C\;,
$$
$$
k'_0\; {\bf l}:\qquad s D=sD+C \alpha\;,\quad k_0\; {\bf
l}':\qquad Ds+\alpha=D\;,
$$
$$
l'_0\; {\bf l}:\qquad tD +B=Dt+C\beta\;,\quad l_0\; {\bf
l}':\qquad Dt+\beta= D^2 +CB\;,
$$
$$
0\; {\bf k}:\qquad t C+A=D\;,\quad l_0\; {\bf k}':\qquad t C =DC
+CA\;,
$$
$$
{\bf k}'\times {\bf k} :\qquad C^2=C\;,\quad {\bf l}'\times {\bf
l} :\qquad D^2+B=D^2+CB\;,
$$
$$
{\bf k}'\times {\bf l} :\qquad CD=CD +CA\;,\quad {\bf l'}\times
{\bf k} :\qquad DC +A=D\;.
$$

The equation
\begin{eqnarray}m_0''=t l_0''+s k_0''\label{B9.7a}\end{eqnarray}
gives
$$(t l_0'+s k_0') (t l_0+s k_0)+(D {\bf l}'
+ C {\bf k}')\; (D {\bf l}+C {\bf k}) + l'_0\; (\alpha k_0+\beta
l_0)+{\bf l}' (A {\bf k} + B {\bf l})
$$
$$
= t\; [\; l_0'\; k_0+{\bf l}'\; {\bf k} +
 (t l_0'+s k_0') l_0+(D {\bf l}'+C {\bf k}')\; {\bf l}\;]
 $$
 $$
 + s\; [\; k_0'\; k_0+{\bf k}'\cdot{\bf k}+(\alpha k_0'+\beta l_0')\; l_0 +
(A {\bf k}'+B {\bf l}')\cdot {\bf l}\; ]\; ;
$$
whence  it follows
$$
k'_0k_0:\qquad s^2=s\;, \quad l'_0l_0:\qquad t^2+\beta
=t^2+s\beta\;,
$$
$$
k'_0l_0:\qquad st=ts+s \alpha\;,\quad l'_0k_0:\qquad
ts+\alpha=t\;,
$$
$$
{\bf k}'\cdot {\bf k}:\qquad C^2=s\;,\quad {\bf l}' \cdot{\bf
l}:\qquad D^2+ B=t D +s B\;,
$$
$$
{\bf k}' \cdot{\bf l}:\qquad CD =tC +s A\;,\quad {\bf k}\cdot
{\bf l}':\qquad DC +A=t\;.
$$

Let us collect the equations together:
\begin{eqnarray} \alpha C=0\;,\qquad AC=0\;,\qquad AB=0\;,\nonumber\\
    (s-1) B=0\;,\qquad\beta C=B\;,\qquad (C-1) B=0\;,\nonumber\\
    \alpha+A (s-1)=0\;,\qquad (A -D)\alpha+B (s-1)=0\;,\nonumber\\
    (A -D)\beta+B t=0\;,\qquad A B+B(D-t)=0\;,\nonumber\\
  \beta +At=A^2+BC\;,\qquad B+AD=A^2+BC\;,\label{B9.8a}\end{eqnarray}
\begin{eqnarray}\alpha s=0\;,\qquad \alpha\beta=0\;,\nonumber\\
    (s -1)\beta=0\;,\qquad BC=\beta\;,\nonumber\\
    \alpha t=\alpha^2+\beta(s-1)\;,\qquad A+AC=\alpha\;,\nonumber\\
    (B-\beta) D=\alpha B\;,\qquad B +AD=\alpha A+\beta C\;,\label{B9.8b}\end{eqnarray}
\begin{eqnarray} C \alpha=0\;,\qquad CA= 0\,,\nonumber\\
    (s-1) C=0\;,\qquad B=C\beta\;,\qquad C(C-1) =0\;,\qquad (C -1) B=0\;,\nonumber\\
    D(s -1)+\alpha=0\;,\qquad Dt+\beta= D^2 +CB\;,\nonumber\\
    t C+A=D\;,\qquad (D-t)C +CA=0\;,\qquad DC +A=D\;;\label{B9.8c}\end{eqnarray}
\begin{eqnarray} 0=s \alpha\;,\qquad  s(s-1)=0\;,\qquad (s-1)\beta=0\;,\nonumber\\
    t(s -1)+\alpha=0\;,\qquad C^2=s\;,\nonumber\\
    D(D- t)=(s -1) B\;,\qquad C(D -C)=s A\;,\qquad DC +A=t\;.\label{B9.8d}\end{eqnarray}

For $s$, only two values are possible: $s=0,\,1$. For $s=0$, the
system has only one solution, namely

\vspace{3mm}

\underline{solution $(KL-1)$},
\begin{eqnarray}\alpha=A\;,\qquad B=\beta=0\;,\qquad C= s=0\;,\qquad D=t=A\;,\nonumber\\
 {\bf n}=A {\bf k}\;,\qquad n_0=A k_0\;,\qquad{\bf m}=A {\bf l}\;,\qquad m_0=A l_0\;,
\label{B9.9a}\end{eqnarray}
$$
G=\left(\ba{cc}K&AK\\L&A L\ea\right), \;\; G'G=
\left(\ba{cc}(K'K+AK'L)&A(K'K+AK'L)\\(L'K+ AL'L)&A (L'K+ AL'L)
\ea\right).
$$

For $s=1$, the system also has only one solution:

\vspace{3mm}

\underline{solution $(KL-2)$},
\begin{eqnarray}C= s= +1\;,\qquad A=\alpha=0\;,\qquad D=t=+1\;,\qquad B=\beta =1\;,\nonumber\\
 {\bf n}= {\bf l}\;,\qquad n_0= l_0\;,\qquad {\bf m}= {\bf l}+ {\bf k}\;,\qquad  m_0= l_0+ k_0\;,
\label{B9.10a}\end{eqnarray}
$$
G=\left(\ba{cc} K&L\\ L&K+L \ea\right),
$$
$$
G'G= \left(\ba{cc}(K'K +L'L)&(L'K+ K'L+L'L)\\
(L'K+ K'L+L'L)&[(K'K+L'L) +(L'K+K'L+L'L)]\ea\right).
$$

\section{Two Independent Vectors: Variant I(n,m)}
Let us examine {\bf variant II(n,m)}:
\begin{eqnarray} {\bf l}=A {\bf m}+B {\bf n}\;,\qquad n_0=\alpha m_0+\beta n_0\;,\nonumber\\
{\bf k}=D {\bf n}+C {\bf m}\;,\qquad k_0=t n_0+s
m_0\;.\label{B10.1}\end{eqnarray}

Here we get only

\vspace{3mm}

\underline{solution $(NM-1)$},
\begin{eqnarray}G=\left(\ba{cc}AN&N\\AM&M\ea\right) ;\label{B10.2}\end{eqnarray}

 \underline{solution $(NM-2)$},
\begin{eqnarray}G= \left(\ba{cc}M+N&N\\N&M\ea\right).\label{B10.3}\end{eqnarray}
\section{Two Independent Vectors: Variant II(k,n)}
Let us examine {\bf variant II(k,n)}:
\begin{eqnarray} {\bf l}=(A {\bf k}+B {\bf n})\;, \qquad l_0=(\alpha k_0+\beta n_0)\;,\nonumber\\
    {\bf m}=(D {\bf n}+C {\bf k})\;,\qquad m_0=(t n_0+s k_0)\;.\label{B11.1}\end{eqnarray}
The multiplication law takes the form
\begin{eqnarray}k_0''=k_0'\; k_0+{\bf k}'\cdot {\bf k}+n'_0\; (\alpha
k_0+\beta n_0)+{\bf n}'\; (A {\bf k}+B {\bf n})\;,\nonumber\\
m_0''=(t n_0'+s k_0')\; (t n_0+s k_0) +
(D {\bf n}'+C {\bf k}')\; (D {\bf n}+C {\bf k})\nonumber\\
 + (\alpha k_0'+\beta n_0') n_0+(A {\bf k}'+B {\bf n}')\; {\bf n}\;,\nonumber\\
n_0''=k_0'\; n_0+{\bf k}'\cdot {\bf n}
+n'_0\; (t n_0+s k_0)+{\bf n}'\; (D {\bf n}+C {\bf k})\;,\nonumber\\
l_0''=(\alpha k_0'+\beta n_0')\; k_0+(A {\bf k}'+B {\bf n}')\; {\bf k}\nonumber\\
 + (t n_0'+s k_0') (\alpha k_0+\beta n_0) + (D {\bf n}'+C {\bf k}')\; (A {\bf k}+B {\bf n})\;,
\nonumber \end{eqnarray}
\begin{eqnarray}{\bf k}''=k'_0\; {\bf k}+{\bf k}'\; k_0+i\; {\bf
k}'\times {\bf k}+n_0'\; (A {\bf k}\nonumber\\
+B {\bf n})+{\bf n}'\; (\alpha k_0+\beta n_0)+i\; {\bf n}'\times (A {\bf k}+B {\bf n})\;,\nonumber\\
{\bf m}''=(t n_0'+s k_0') (D {\bf n}+C {\bf k})+ (D {\bf n}'+C {\bf k} ')\; (t n_0+s k_0)\nonumber\\
+i\;(D {\bf n}'+C {\bf k} ')\times (D {\bf n}+C {\bf k})\nonumber\\
 + (\alpha k_0'+\beta n_0')\; {\bf n}+(A {\bf k}'+B {\bf n}')\; n_0+i\;
 (A {\bf k}'+B {\bf n}')\times {\bf n}\;,\nonumber\\
{\bf n}''=k'_0\; {\bf n}+{\bf k}'\; n_0+i\; {\bf
k}'\times {\bf n}+n_0'\; (D {\bf n}+C {\bf k})\nonumber\\
+{\bf n}'\; (t n_0+s k_0)+i\; {\bf n}'\times (D {\bf n}+C {\bf k})\;,\nonumber\\
{\bf l}''=(\alpha k_0'+\beta n_0')\; {\bf k}+(A {\bf k}'+B {\bf
n}')\; k_0+i\;
 (A {\bf k}'+B {\bf n}')\times {\bf k} \nonumber\\
 + (t n_0'+s k_0') (A {\bf k}+B {\bf n})+(D {\bf n}'+C {\bf k}') (\alpha k_0+\beta n_0)\nonumber\\
 + i\; (D {\bf n}'+C {\bf k}')\times (A {\bf k}+B {\bf n})\;.\label{B11.2} \end{eqnarray}

The equation
\begin{eqnarray} {\bf l}''=(A {\bf k}''+B {\bf n}'')\label{B11.3a} \end{eqnarray}
is equivalent to
\begin{eqnarray}(\alpha k_0'+\beta n_0') {\bf k}+(A {\bf k}'+B
{\bf n}') k_0+i (A {\bf k}'+B {\bf n}')\times {\bf k}\nonumber\\
 + (t n_0'+s k_0') (A {\bf k}+B {\bf n})+(D {\bf n}'+C {\bf k}') (\alpha k_0+\beta n_0) +
 i (D {\bf n}'+C {\bf k}')\times (A {\bf k}+B {\bf n})\nonumber\\
=A [ k'_0 {\bf k}+{\bf k}' k_0+i {\bf k}' \times {\bf k}+n_0' (A
{\bf k}+B {\bf n})+{\bf n}' (\alpha k_0+\beta n_0)+i {\bf
n}'\times (A {\bf k} +
B {\bf n}) ]\nonumber\\
+ B [ k'_0\; {\bf n}+{\bf k}' n_0+i {\bf k}' \times {\bf n}+n_0'
(D {\bf n}+C {\bf k})+{\bf n}' (t n_0+s k_0)+i {\bf n}'\times (D
{\bf n} + C {\bf k}) ] ;\nonumber\end{eqnarray}
this leads to
$$
k'_0 {\bf k}:\qquad \alpha+s A=A\;,\quad k_0 {\bf k}':\qquad
A+\alpha C=A\;,
$$
$$
n_0 {\bf n}':\qquad\beta D=\beta A+t B\;,\quad n_0' {\bf n}:\qquad
t B=AB+B D\;,
$$
$$
k'_0 {\bf n}:\qquad s B =B\;,\quad k_0 {\bf n}':\qquad B+\alpha
D=\alpha A+s B\;,
$$
$$
n'_0 {\bf k}:\qquad\beta +tA=A^2+BC\;,\quad n_0 {\bf
k}':\qquad\beta C=B\;,
$$
$$
{\bf k}'\times {\bf k} :\qquad A+AC=A\;,\quad {\bf n}'\times {\bf
n} :\qquad BD=AB +BD\;,
$$
$$
{\bf n}'\times {\bf k} :\qquad B+AD=A^2+BC\;,\quad {\bf k}'\times
{\bf n} :\qquad CB=B\;.
$$

The equation
\begin{eqnarray}{\bf m}''=(D {\bf n}''+C {\bf k} '')\label{B11.4a}\end{eqnarray}
takes the form
\begin{eqnarray} (t n_0'+s k_0') (D {\bf n}+C {\bf k})+(D {\bf n}'+C {\bf k} ')
(t n_0+s k_0)+i (D {\bf n}'+C {\bf k} ')
\times (D {\bf n}+C {\bf k})\nonumber\\
 + (\alpha k_0'+\beta n_0') {\bf n}+(A {\bf k}'+B {\bf n}') n_0+i
(A {\bf k}'+B {\bf n}')\times {\bf n}\nonumber\\
= D [ k'_0 {\bf n}+{\bf k}' n_0+i {\bf k}' \times {\bf n}+n_0' (D
{\bf n}+C {\bf k})+{\bf
n}' (t n_0+s k_0)+i {\bf n}'\times (D {\bf n} +C {\bf k}) ]\nonumber\\
+ C [ k'_0 {\bf k}+{\bf k}' k_0+i {\bf k}' \times {\bf k}+n_0' (A
{\bf k}+B {\bf n})+{\bf n}' (\alpha k_0+\beta n_0)+i {\bf
n}'\times (A {\bf k} + B {\bf n}) ] ;\nonumber\end{eqnarray}
whence it follows
$$
k'_0 {\bf k}:\qquad sC=C\;,\quad
 k_0 {\bf k}':\qquad sC=C\;,
 $$
 $$
    n_0 {\bf n}':\qquad tD +B=t D+\beta C\;,\quad
     n_0' {\bf n}:\qquad tD+\beta =D^2+CB\;,
     $$
     $$
     k'_0 {\bf n}:\qquad sD+\alpha=D\;,\quad
     k_0 {\bf n}':\qquad sD=sD+\alpha C\;,
     $$
     $$
     n'_0 {\bf k}:\qquad t C=DC+CA\;,\quad
      n_0 {\bf k}':\qquad tC +A=D\;,
      $$
      $$
          {\bf k}'\times {\bf k} :\qquad C^2=C\;,\quad
          {\bf n}'\times {\bf n} :\qquad D^2+B=D^2+CB\;,
          $$
          $$
          {\bf n}'\times {\bf k} :\qquad DC=DC+CA\;,\quad
          {\bf k}'\times {\bf n} :\qquad CD +A=D\;.
          $$

The equation
\begin{eqnarray}l_0''=(\alpha k_0''+\beta n_0'')\label{B11.5a}\end{eqnarray}
gives
\begin{eqnarray} (\alpha k_0'+\beta n_0')\; k_0+(A {\bf k}'+B {\bf n}')\; {\bf k}\nonumber\\
    + (t n_0'+s k_0') (\alpha k_0+\beta n_0)+(D {\bf n}'+C {\bf k}')\; (A {\bf k}+B {\bf n})\nonumber\\
    =\alpha\; [\; k_0'\; k_0+{\bf k}'\cdot {\bf k}+n'_0\;(\alpha k_0+\beta n_0)+{\bf n}'\; (A {\bf k}+B
    {\bf n})\;]\nonumber\\
    +\beta\; [\; (t n_0'+s k_0')\; (t n_0+s k_0)+ (D {\bf n}'+C {\bf k}')\; (D {\bf n}+
    C {\bf k})\;]\;,\nonumber\end{eqnarray}
which leads to
$$
k_0'k_0:\qquad \alpha+\alpha s=\alpha+\beta s^2\;,\quad
    n_0'n_0:\qquad t\beta=\alpha\beta+\beta t^2\;,
    $$
    $$
    k_0'n_0:\qquad s\beta=\beta s t\;,\quad
    k_0n_0':\qquad\beta+t \alpha=\alpha^2+\beta ts\;,
    $$
    $$
    {\bf k} ' \cdot{\bf k} :\qquad A+CA=\alpha+\beta C^2\;,\quad
    {\bf n} '\cdot {\bf n} :\qquad DB=\alpha B+\beta D^2\;,
    $$
    $$
        {\bf k} ' \cdot{\bf n} :\qquad CB=\beta CD\;,\quad
        {\bf k} \cdot{\bf n}' :\qquad B +AD=\alpha A+\beta CD\;.
        $$

The equation
\begin{eqnarray}m_0 ''= t n_0''+s k_0''\label{B11.6a}\end{eqnarray}
is equivalent to
\begin{eqnarray} (t n_0'+s k_0')\; (t n_0+s k_0) +(D {\bf n}'+C {\bf k}')\; (D {\bf n}+C {\bf k})\nonumber\\
    + (\alpha k_0'+\beta n_0') n_0+(A {\bf k}'+B {\bf n}')\; {\bf n}\nonumber\\
    = t\; [\; k_0'\; n_0+{\bf k}'\cdot {\bf n}+n'_0\; (t n_0+s k_0)+{\bf n}'\; (D {\bf n}+C {\bf k})\nonumber\\
    + s [\;k_0'\; k_0+{\bf k}'\cdot {\bf k}+n'_0\; (\alpha k_0 +\beta n_0)+{\bf n}'\; (A {\bf k}+B {\bf n})\; ]
    \; ;\nonumber\end{eqnarray}
whence it follows
$$
k_0'k_0:\qquad s^2 =s\;,\quad
 n_0'n_0:\qquad t^2+\beta=t^2+s\beta\;,
 $$
 $$
     k_0'n_0:\qquad st+\alpha=t\;,\quad
     k_0n_0':\qquad ts =st+\alpha s\;,
     $$
     $$
         {\bf k} '\cdot {\bf k} :\qquad C^2=s\;,\quad
          {\bf n} ' \cdot{\bf n} :\qquad D^2+B=tD+s B\;,
          $$
          $$
              {\bf k} '\cdot {\bf n} :\qquad CD+A=t\;,\quad
               {\bf n}' \cdot{\bf k} :\qquad DC=t C+sA\;.
               $$

We collect results together, and yield
$$
\alpha+s A=A\;,\; \alpha C=0\;,\;\beta D=\beta A+t B\;,\;
    t B=AB+B D\;,\;  s B =B\;,
    $$
    $$
     B+\alpha D=\alpha A+s B\;,\;\beta +tA=A^2+BC\;,\;
  \beta C=B\;,\; AC=0\;,\; 0=AB\;,
  $$
  $$
   B+AD=A^2+BC\;,\;  CB=B\; ;\;   sC=C\;,\;     sC=C\;,
   $$
   $$
     B=\beta C\;,\;  tD+\beta =D^2+CB\;,\;     sD+\alpha=D\;,\; 0=\alpha C\;,
     $$
     $$
      t C=DC+CA,\; tC +A=D,\;
    C(C-1)=0,\;
    $$
    $$
    B=CB,\;  0= CA,\;
    CD +A=D , \;\alpha s=\beta s^2,
    $$
$$
t\beta=\alpha\beta+\beta t^2,\; s\beta=\beta s t,\;
  \beta+t \alpha=\alpha^2+\beta ts,\; A+CA=\alpha+\beta C^2,
  $$
  $$
    DB=\alpha B+\beta D^2,\; CB=\beta CD,\; B +AD=\alpha A+\beta CD ;
$$
$$
 s(s-1) =0\;,\; \beta=s\beta\;,\;
    st+\alpha=t\;,\;  0 =\alpha s\;, \; C^2=s\;,
$$
$$
    D^2+B=tD+s B\;,\;  CD+A=t\;,\;  DC=t C+sA\;.
    $$

We see that for $s,C$ there are possible only the following values
\begin{eqnarray}C=s=0\;,\qquad C=s=1\;.\nonumber\end{eqnarray}

If $C=s=0$, we get only one

\vspace{3mm}

\underline{solution $(KN-1)$},
\begin{eqnarray} \alpha=A\;,\quad B=\beta =0\;,\quad C=s=0\;,\quad D=t=A\;,\nonumber\\
 {\bf l}=A {\bf k}\;,\quad l_0=A k_0\;,\quad {\bf m}=A {\bf n}\;,\;
 m_0=A n_0\;,\label{B11.7b}\end{eqnarray}
 $$
 G=\left(\ba{cc}K&N\\AK&AN\ea\right),\; G'G=
 \left(\ba{cc}(K'K +A N'K)&(K'N +A N'N)\\A(K'K +A N'K)&A(K'N +A N'N)
\ea\right).
$$

In the case $C=s=+1$, we also get only one

\vspace{3mm}

\underline{solution $(KN-2)$},
\begin{eqnarray}
A= \alpha =0\;,\quad B=\beta =0\;,\quad C=s=+1\;,\quad D= t=0\;,
\nonumber\\
{\bf l}=0\;,\quad l_0=0\;,\quad {\bf m}= {\bf k}\;,\quad m_0=
k_0\;,\label{B11.8a}
\end{eqnarray}
$$
G=\left(\ba{cc}K&N\\0&K\ea\right),\;G'G= \left(\ba{cc}K'K&(K'N+
N'K)\\0&K'K \ea\right).
$$

\section{Two Independent Vectors: Variant II(m,l)}
Let us examine {\bf variant II(m,l)}:
\begin{eqnarray} {\bf n}=A {\bf m}+B {\bf l}\;,\qquad l_0=\alpha m_0+\beta l_0\;,\; \nonumber\\
    {\bf k}=D {\bf l}+C {\bf m}\;,\qquad k_0=t l_0+s m_0\;.\;\; \label{B12.1}\end{eqnarray}
Here there arise only

\vspace{3mm}

\underline{solution $(ML-1)$},
\begin{eqnarray}G=S G'S^{-1}=\left(\ba{cc}AN'&AK'\\N'&K '\ea\right) \equiv
    \left(\ba{cc}AL&A M\\L&M\ea\right) ;\label{B12.2}\end{eqnarray}

 \underline{solution $(ML-2)$},
\begin{eqnarray}G=SG' S^{-1}=\left(\ba{cc}K'&0\\N'&K'\ea\right) \equiv
    \left(\ba{cc}M&0\\L&M\ea\right).\label{B12.3}\end{eqnarray}
\section{Three Independent Vectors: Variants  I(k,m,n) and I(m,k,l)}
Now consider the  case of three independent vectors
\begin{eqnarray}{\bf l}=A {\bf k}+B {\bf m}+C {\bf n}\;,\qquad l_0=\alpha k_0+
  \beta m_0+s n_0\;.\label{B13.1}\end{eqnarray}

First, let us examine the simpler second restriction in
(\ref{B13.1}). The formulas for
    $k''_0,\;m_0'',\;n''_0,\;l''_0$ take the form
\begin{eqnarray} k_0''=k_0' k_0+{\bf k}' \cdot{\bf k} +n'_0 (\alpha k_0 +
  \beta m_0+s n_0)+{\bf n}' (A {\bf k}+B {\bf m}+C {\bf n})\;,\qquad\nonumber\\
    m_0''=m_0' m_0+{\bf m}'\cdot {\bf m} +(\alpha k'_0+\beta m'_0+s n'_0) n_0 +
    (A {\bf k}'+B {\bf m}'+C {\bf n}') {\bf n}\;, \nonumber\end{eqnarray}
\begin{eqnarray}n_0''=k_0'\; n_0+{\bf k}'\cdot {\bf n}+n'_0\; m_0+{\bf n}'\cdot {\bf m}\;,\nonumber\\
    l_0''=(\alpha k'_0+\beta m'_0+s n'_0) k_0+(A {\bf k}+B {\bf m}+C {\bf n}) {\bf k}\nonumber\\
    +m'_0 (\alpha k_0+\beta m_0+s n_0)+{\bf m}'  (A {\bf k}+B {\bf m}+C {\bf n})\;. \nonumber\end{eqnarray}

We set
\begin{eqnarray}l''_0=\alpha k''_0+\beta m''_0+s n''_0\;\label{B13.2a}\end{eqnarray}
or, equivalently,
\begin{eqnarray}(\alpha k'_0+\beta m'_0+s n'_0) k_0+(A {\bf k}' +B {\bf m}'+C {\bf n}') {\bf k}\nonumber\\
+m'_0 (\alpha k_0+\beta m_0+s n_0)+{\bf m}' (A {\bf k}+B {\bf m}+C {\bf n})\nonumber\\
= \alpha\; [\;k_0' k_0+{\bf k}'\cdot {\bf k} +n'_0 (\alpha
k_0+\beta
m_0+s n_0)+{\bf n}' (A {\bf k}+B {\bf m}+C {\bf n})\; ]\nonumber\\
+\beta[\; m_0' m_0+{\bf m}'\cdot {\bf m}+(\alpha k'_0+\beta m'_0+s
n'_0) n_0 +
(A {\bf k}'+B {\bf m}'+C {\bf n}') {\bf n}\; ]\nonumber\\
+ s\; [\;k_0'\; n_0+{\bf k}'\cdot {\bf n}+n'_0\; m_0+{\bf n}'\cdot
{\bf m}\; ]\;,\nonumber\end{eqnarray}
which results in
\begin{eqnarray}k'_0k_0:\; \alpha=\alpha\;,\qquad\qquad m'_0m_0:\;\beta=\beta ;,\nonumber\\
m_0'k_0 :\;\beta+\alpha =0\;,\qquad\qquad n'_0n_0:\; 0=(\alpha+\beta) s\;,\nonumber\\
n'_0 k_0:\; s=\alpha^2\;,\qquad\qquad m'_0 n_0:\; s=\beta^2\;,\nonumber\\
n'_0m_0:\; 0=\beta (\alpha+s)\;,\qquad k'_0 n_0:\; 0=\beta (\alpha +s)\;,\nonumber\\
{\bf k}'\cdot{\bf k} :\; A =\alpha\;,\qquad\qquad {\bf m}'\cdot{\bf m}:\; B =\beta\;,\nonumber\\
{\bf m}'\cdot{\bf k} :\; B+A=0\;,\qquad\qquad {\bf n}'\cdot{\bf n} :\; 0=(\alpha+\beta) C\;,\nonumber\\
{\bf n}'\cdot{\bf k} :\; C=\alpha A\qquad\qquad {\bf m}'\cdot{\bf n}:\; C=\beta B\;,\nonumber\\
{\bf n}'\cdot{\bf m} :\; 0=\alpha B+\beta s\;,\qquad {\bf
k}'\cdot{\bf n} :\; 0=\beta A+\beta
s\;.\label{B13.2b}\end{eqnarray}

Let us write down only the independent equations:
\begin{eqnarray}\alpha=A,\qquad B=\beta=-A,\qquad C=s=A^2,
    \qquad\;\;\; A(A+s)=0\;.\label{B13.2c}\end{eqnarray}
The equations (\ref{B13.2c}) allow two types of solutions:
\begin{eqnarray}(I),\qquad A=\alpha=0\;,\qquad B=\beta =0\;,\; C= s=0\; ;\label{B13.3a}\\
    (II),\qquad A= \alpha=-1,\qquad B=\beta=+1\;,\qquad C= s=+1\;.\label{B13.3b} \end{eqnarray}

It remains to verify that these solutions satisfy the equation
(see (\ref{B13.1}))
\begin{eqnarray}{\bf l}''=A {\bf k} ''+ B {\bf m}''+C {\bf n}''\;.\label{B13.4a}\end{eqnarray}
Evidently, we should check only the type $II$ -- it takes the form
\begin{eqnarray}{\bf l}'' =-{\bf k} ''+ {\bf m}''+{\bf n}''\;\nonumber\end{eqnarray}
and further reduces to
\begin{eqnarray}(- k'_0+m'_0+n'_0)\;{\bf k}+(- {\bf k}'+{\bf m}'+{\bf n}')\; k_0\nonumber\\
 + i\; (- {\bf k}'+{\bf m}'+{\bf n}')\times {\bf k} +m_0'\; (- {\bf k}+{\bf m}+{\bf n})\nonumber\\
 + {\bf m}' (- k_0+m_0+n_0)+i\; {\bf m}'\times (-{\bf k}+{\bf m}+{\bf n})\nonumber\\
=- [\; k'_0\; {\bf k}+{\bf k}'\; k_0+i\; {\bf k}'\times {\bf k}+n_0'\; (- {\bf k}+{\bf m}+{\bf n})\nonumber\\
+ {\bf n}'\;(- k_0+m_0+n_0)+i\;{\bf n}'\times (- {\bf k}+{\bf m}+{\bf n})\; ]\nonumber\\
 + [\; m'_0\; {\bf m}+{\bf m}'\; m_0+i\; {\bf m}'\times {\bf m}+(- k_0'+m_0'+n_0')\; {\bf n} +\nonumber\\
 (- {\bf k}'+{\bf m}'+{\bf n}')\; n_0+i\;(- {\bf k}'+{\bf m}'+{\bf n}')\times {\bf n}\; ]\nonumber\\
 + [\; k'_0\; {\bf n}+{\bf k}'\; n_0+i\; {\bf k}'\times {\bf n}+n_0'\; {\bf m}+{\bf n}'\; m_0+i\;
{\bf n}'\times {\bf m}\; ]\;,\nonumber\end{eqnarray}
which turns to be an identity.

Thus, we have constructed two solutions:

\vspace{3mm}

\underline{ solution $(KMN-1)$},
\begin{eqnarray} G=\left(\ba{cc} K&N\\ 0&M \ea\right),\;
    G'G=
    \left(\ba{cc} K' K&(K'N +N'M)\\ 0&M'M \ea\right) ;\label{B13.5}\end{eqnarray}

\underline{solution $(KMN-2)$},
\begin{eqnarray} {\bf l} =-{\bf k}+{\bf m}+{\bf n}\;,\;l_0 =-k_0+m_0+n_0\;,\nonumber\\
    G=\left(\ba{cc} K&N\\ -K +M+N&M \ea\right) ;\label{B13.6a}\end{eqnarray}
These are sets of degenerate matrices of rank 2. It remains to
verify the multiplication law for the case (\ref{B13.6a}):
\begin{eqnarray}\left(\ba{cc} K'&N'\\ -K' +M'+N '&M'\ea\right)
    \left(\ba{cc} K&N\\ -K +M+N&M \ea\right) \qquad \quad
         \nonumber\\
    =\left(\ba{cc}(K'K-N'K +N' M+N'N)&K'N+ N'M\\
    (-K'K +N ' K+M'M +M' N)&(-K'N +M'N+N ' N+M'M)
    \ea\right)\;.\nonumber\\\label{B13.6b}\end{eqnarray}
With the notations
$$K''=K'K-N'K +N' M+N'N\;,
$$
$$
    M''= -K'N +M'N+N ' N+M'M\;,\quad N''=K'N+ N'M\;,
    $$
    we obtain
$$
-K''+M''+N''=-K'K+N'K -N' M\; ,
$$
$$
    - N'N -K'N +M'N+N ' N+M'M+K'N+ N'M = -K'K+N'K +M'N+M'M\; ;
    $$
    Thus, the set (\ref{B13.6a}) has semi-group structure
\begin{eqnarray}G'G=\left(\ba{cc} K''&N''\\ -K'' +M''+N ''&M'
    \ea\right).\label{B13.6c}\end{eqnarray}

For the variant  {\bf II(k,m,l)} one has the similar solutions:

\vspace{3mm}

\underline{ solution $(KML-1)$},
\begin{eqnarray}\qquad G=
    \left(\ba{cc} K&0\\ L&M \ea\right) ;\label{B13.7a}\end{eqnarray}

\underline{ solution $(KML-2)$},
\begin{eqnarray} G= \left(
    \ba{cc}K&-M +K+L\\ L&M \ea\right).\label{B13.7b}\end{eqnarray}
\section{Three Independent Vectors: Variants I(n,l,k) and I(n,l,m)}
Let us consider the variant {\bf  $I(n,l,k)$}
%
\begin{eqnarray}{\bf m}=A {\bf n}+B {\bf l}+C {\bf k}\;,\qquad m_0=
    \alpha n_0+\beta l_0+s k_0\;.\label{B14.1}\end{eqnarray}

First, we turn to the simplest second restriction in
(\ref{B14.1}); this gives
\begin{eqnarray}k_0''=k_0'\; k_0+{\bf k}'\cdot {\bf k}+n'_0\; l_0
+ {\bf n}'\cdot {\bf l}\;,\qquad\qquad\nonumber\\
m_0''=(\alpha n'_0+\beta l'_0+s k'_0) (\alpha n_0
+\beta l_0+s k_0)\qquad\qquad\nonumber\\
+ (A {\bf n}'+B {\bf l}'+C {\bf k}') (A {\bf n}+B {\bf l}+C {\bf
k})
+l'_0\; n_0+{\bf l}'\cdot {\bf n}\;,\qquad\qquad\nonumber\\
n_0''=k_0'\; n_0+{\bf k}'\cdot {\bf n}+n'_0 (\alpha n_0+\beta
l_0+s k_0)+{\bf n}'\;
(A {\bf n}+B {\bf l}+C {\bf k})\;,\nonumber\\
l_0''=l_0'\; k_0+{\bf l}'\cdot {\bf k} +(\alpha n'_0+\beta l'_0+s
k'_0) l_0+(A {\bf n}'+B {\bf l}'+C {\bf k}'){\bf l}\;.
\label{B14.2a}\end{eqnarray}

By requiring
\begin{eqnarray}m_0''=\alpha n''_0+\beta l''_0+s k''_0\;,\nonumber\end{eqnarray}
or
\begin{eqnarray}(\alpha n'_0+\beta l'_0+s k'_0) (\alpha n_0+\beta l_0+s k_0)\nonumber\\
 + (A {\bf n}'+B {\bf l}'+C {\bf k}') (A {\bf n}+B {\bf l}+C {\bf k})
+l'_0\; n_0+{\bf l}'\cdot {\bf n}\nonumber\\
= \alpha\; [\; k_0'\; n_0+{\bf k}'\cdot {\bf n} +n'_0 (\alpha
n_0+\beta l_0+s k_0)+{\bf n}'\;
 (A {\bf n}+B {\bf l}+C {\bf k})\; ]\nonumber\\
 +\beta\; [\; l_0'\; k_0+{\bf l}'\cdot {\bf k}
+(\alpha n'_0+\beta l'_0+s k'_0) l_0+(A {\bf n}'+B {\bf l}'+C {\bf
k}')
 {\bf l}\; ]\nonumber\\
 + s\; [\; k_0'\; k_0+{\bf k}'\cdot {\bf k} +n'_0\; l_0+{\bf n}'\cdot {\bf l}\; ]\;,\nonumber \end{eqnarray}
we obtain the equations
\begin{eqnarray}n_0'n_0:\; \alpha^2=\alpha^2\;,\qquad\qquad
l'_0l_0:\;\beta^2=\beta^2\;,\qquad k'_0k_0:\;
s^2=s\;,\nonumber\\
n'_0l_0:\; \alpha\beta=\alpha\beta+\alpha\beta+s\;
,\qquad n_0'k_0:\; \alpha s=\alpha s\;,\nonumber\\
l'_0n_0:\;\beta \alpha +1=0\;,\qquad l'_0 k_0:\;
\beta s=\beta\,,\nonumber\\
k'_0n_0:\; s \alpha =\alpha\;,\qquad k_0'l_0:\; s \beta=\beta
s\;,\label{B14.2b}\end{eqnarray}
\begin{eqnarray}{\bf n}'\cdot {\bf n}:\; A^2 =\alpha A,\qquad {\bf l}'\cdot{\bf l}:\;
B^2=\beta B\;,\qquad {\bf k}'\cdot{\bf k}:\; C^2 =s\;,\nonumber\\
{\bf n}' \cdot{\bf l}:\; AB=\alpha B+\beta A+s\;,\qquad {\bf
n}' \cdot{\bf k}:\; AC=\alpha C\;,\nonumber\\
{\bf l}' \cdot{\bf n}:\; BA +1 =0\;,\qquad {\bf l}'\cdot {\bf k}:\; BC=\beta\;,\nonumber\\
{\bf k}' \cdot{\bf n}:\; CA=\alpha\;,\qquad
 {\bf k}'\cdot{\bf l}:\; CB=\beta C\;.\label{B14.2c}\end{eqnarray}
Let us write down only the independent equations:
\begin{eqnarray} \alpha\beta=-1,\qquad s=-\alpha\beta=+1,\qquad\qquad C^2= s=1\;,\nonumber\\
    A(A-\alpha)=0,\qquad C(A-\alpha)=0,\qquad B(B-\beta)=0,\nonumber\\
    AB=\alpha B+\beta A+1\;,\qquad AC=\alpha C\;,\nonumber\\BA +1 =0\;,\quad BC=\beta\;,\quad
    CA=\alpha\;,\quad CB=\beta C\; ;
    \nonumber\end{eqnarray}
they provide to us only with \vspace{3mm}

\underline{solution $(NLK-1)$},
\begin{eqnarray}\alpha=A,\qquad B=\beta=-{1\over A}\;,\qquad
    \qquad C= s=+1,\nonumber\\{\bf m} ={\bf k}+(A\; {\bf n}-A^{-1}\; {\bf l})\;,
    \qquad m_0=k_0+(A\; n_0-A^{-1} l_0)\; ;\label{B14.3a}\end{eqnarray}
where the corresponding matrices have the structure
\begin{eqnarray}G=\left(\ba{cc}K&N\\L&(K +A N -A^{-1} L)
    \ea\right).\label{B14.3b}\end{eqnarray}
This is a set of degenerate matrices of rank 2. We verify now the
multiplication law:
\begin{eqnarray}G'G= \left(\ba{cc}K'&N'\\L '&(K ' +A N' -A^{-1} L')
\ea\right) \left(\ba{cc}K&N\\L&(K +A N -A^{-1} L)\ea\right)\nonumber\\
 \left(\ba{cc}(K'K+N'L)&K'N+N'(K +A N -A^{-1} L)\\ L'K+(K ' +A N' -A^{-1} L') L&L'N +(K ' +A N'
-A^{-1} L')(K +A N -A^{-1}
L)\ea\right).\nonumber\label{B14.4a}\end{eqnarray}
The block (22) in the matrix (\ref{B14.4a}) is
\begin{eqnarray}(22)= L'N +(K ' +A N' -A^{-1} L')\; (K +A N -A^{-1} L)\nonumber\\
    = (K'K-N'L) +A\; (K'N+N'K) - A^{-1}\; (K'L +L'K) +A^2 N'N+A^{-2} L'L\;.\nonumber\end{eqnarray}
With the notations
\begin{eqnarray}K''=K'K+N'L\;,\nonumber\\N''=K'N+N'(K +A N -A^{-1} L)\;,\nonumber\\
L''= L'K+(K ' +A N' -A^{-1} L') L\;,\nonumber\end{eqnarray}
we get
\begin{eqnarray}K''+A N''-A^{-1} L''= K'K+N'L+A\;[\; K'N+N'(K +A N -A^{-1} L)\; ]\nonumber\\
 -A^{-1}\; [\; L'K+(K ' +A N' -A^{-1} L') L\; ]= (K'K- N'L)+A\; (K'N+N'K)\nonumber\\
 -A^{-1}\; (L'K+K ' L) +(A^2 N' N+A^{-2} L' L)\;,\nonumber\end{eqnarray}
and we obtain the claimed relation
\begin{eqnarray}(22)=K'' +A N'' -A^{-1} L''\;.\nonumber\end{eqnarray}
Therefore we have
\begin{eqnarray}G'G= \left(\ba{cc}K''&N''\\L ''&(K '' +A N'' -A^{-1} L'')
    \ea\right)=G''\;.\nonumber\end{eqnarray}

There exists one  solution for the  variant $I(n,l,m)$:

\underline{solution $(NLM-1)$},
\begin{eqnarray} G=
 \left(\ba{cc}(M +A L -A^{-1} N)&N\\L&M
\ea\right)\;.\label{B14.5}\end{eqnarray}
\section{Degenerate Matrices of Rank 3}
Let us consider singular Mueller matrices of rank 3. Given the
explicit form of the
    matrix $ G $, it is easy to understand that there are 16 simple ways to get the
    semigroups of rank 3. For this, it suffices to have vanishing some $i$-line or
    some $j$-column in the original \mbox{4-dimensional} matrix,
\begin{eqnarray}(00)\;,\qquad (01)\;,\qquad (02)\;,\qquad (03)\;,\nonumber\\
(10)\;,\qquad (11)\;,\qquad (12)\;,\qquad (13)\;,\nonumber\\
(20)\;,\qquad (21)\;,\qquad (22)\;,\qquad (23)\;,\nonumber\\
(30)\;,\qquad (31)\;,\qquad (32)\;,\qquad
(33)\;.\label{B15.2}\end{eqnarray}

The compatibility of the law of multiplication with this
constraint is obvious.

All the 16 possibilities are listed below.

\vspace{2mm}

\underline{Variant $ (00) $}
$$G=\left(\ba{cccc}0&0&0&0\\0&2k_0&\;\; 2n_1&2n_0\\
    0&2l_1&\;\; m_0+m_3&m_1-i m_2\\0&2l_0&\;\; m_1+i m_2&m_0-m_3\ea\right),$$
$$k_1 =0,\qquad k_2=0,\qquad k_0=-k_3\;,$$
$$n_0= -n_3,\qquad l_0= -l_3,\qquad +in_2=n_1,\qquad -il_2=l_1\;.$$

\underline{Variant $ (01) $}
$$G=\left(\ba{cccc}0&0&\;\; 0&0\\2k_1&0&\;\; 2 n_1&2n_0\\
    2l_0&0&\;\; m_0+m_3&m_1-i m_2\\2l_1&0&\;\; m_1+i m_2&m_0-m_3\ea\right),$$
$$k_0 =0,\qquad k_3=0,\qquad k_1=ik_2\;,$$
$$l_1= il_2,\qquad l_0= l_3,\qquad n_0=-n_3,\qquad n_1=in_2\;.$$

\underline{Variant $ (02) $}
$$G=\left(\ba{cccc}0&0&\;\; 0&0\\2k_1&2k_0&\;\; 0&2n_0\\
    l_0+l_3&l_1-i l_2&\;\; 0&2m_1\\l_1+i l_2&l_0-l_3&\;\; 0&2m_0\ea\right),$$
$$n_1 =0,\qquad n_2=0,\qquad n_0=-n_3\;,$$
$$m_0= -m_3,\qquad m_1= -im_2,\qquad k_0=-k_3,\qquad k_1=ik_2\;.$$

\underline{Variant $ (03) $}
$$G=\left(\ba{cccc}0&0&\;\; 0&0\\2k_1&2 k_0&\;\; 2n_1&0\\
    l_0+l_3&l_1-i l_2&\;\; 2m_0&0\\l_1+i l_2&l_0-l_3&\;\; 2 m_1&0\ea\right),$$
$$n_0 =0,\qquad n_3=0,\qquad n_1 =i n_2\;,$$
$$m_0= m_3,\qquad m_1= im_2,\qquad k_0=-k_3,\qquad k_1=ik_2\;.$$

\underline{Variant $ (10) $}
$$G=\left(\ba{cccc}0&2k_1&\;\; 2n_0&2n_1\\0&0&\;\; 0&0\\
    0&2l_1&\;\; m_0+m_3&m_1-i m_2\\0&2l_0&\;\; m_1+i m_2&m_0-m_3\ea\right),$$
$$k_0 =0,\qquad k_3=0,\qquad k_1 =-i k_2\;,$$
$$l_0= -l_3,\qquad l_1= -il_2,\qquad n_1=-in_2,\qquad n_0=n_3\;.$$

\underline{Variant $ (11) $}

$$G=\left(\ba{cccc}2k_0&0&\;\; 2n_0&2 n_1\\
    0&0&\;\; 0&0\\2l_0&0&\;\; m_0+m_3&m_1-i m_2\\
    2l_1&0&\;\; m_1+i m_2&m_0-m_3\ea\right),$$
$$k_1 =0,\qquad k_2=0,\qquad k_0=k_3\;,$$
$$l_0= l_3,\qquad l_1= il_2,\qquad n_1 =- in_2,\qquad n_0=n_3\;.$$

\underline{Variant $ (12) $}
$$G=\left(\ba{cccc}2k_0&2k_1&\;\; 0&2n_1\\
    0&0&\;\; 0&0\\l_0+l_3&l_1-i l_2&\;\; 0&2m_1\\
    l_1+i l_2&l_0-l_3&\;\; 0&2 m_0\ea\right),$$
$$n_0 =0,\qquad n_3=0,\qquad n_1 =-i n_2\;,$$
$$m_0= -m_3,\qquad m_1= -im_2,\qquad k_1 =- ik_2,\qquad k_0=k_3\;.$$

\underline{Variant $ (13) $}
$$G=\left(\ba{cccc}2k_0&2 k_1&\;\; 2n_0&0\\
    0&0&\;\; 0&0\\l_0+l_3&l_1-i l_2&\;\; 2m_0&0\\
    l_1+i l_2&l_0-l_3&\;\; 2m_1&0\ea\right),$$
$$n_1 =0,\qquad n_2=0,\qquad n_0=n_3\;,$$
$$m_0= m_3,\qquad m_1= im_2,\qquad k_1 =- ik_2,\qquad k_0=k_3\;.$$

\underline{Variant $ (20) $}
$$G=\left(\ba{cccc}0&2k_1&\;\; n_0+n_3&n_1-i n_2\\
    0&2k_0&\;\; n_1+i n_2&n_0-n_3\\0&0&\;\; 0&0\\
    0&2l_0&\;\; 2m_1&2m_0\ea\right),$$
$$l_1 =0,\qquad l_2=0,\qquad l_0 =- l_3\;,$$
$$m_0= -m_3,\qquad m_1= im_2,\qquad k_1 =- ik_2,\qquad k_0=-k_3\;.$$

\underline{Variant $ (21) $}
$$G=\left(\ba{cccc}2k_0&0&\;\; n_0+n_3&n_1-i n_2\\
    2k_1&0&\;\; n_1+i n_2&n_0-n_3\\0&0&\;\; 0&0\\
    2l_1&0&\;\; 2m_1&2m_0\ea\right),$$
$$l_0 =0,\qquad l_3=0,\qquad l_1 =i l_2\;,$$
$$m_0= -m_3,\qquad m_1= im_2,\qquad k_1=ik_2,\qquad k_0=k_3\;.$$

\underline{Variant $ (22) $}
$$G=\left(\ba{cccc}k_0+k_3&k_1-i k_2&\;\; 0&2n_1\\
    k_1+i k_2&k_0-k_3&\;\; 0&2n_0\\0&0&\;\; 0&0\\
    2l_1&2l_0&\;\; 0&2m_0\ea\right),$$
$$m_1 =0,\qquad m_2=0,\qquad m_0 =-m_3\;,$$
$$n_0= -n_3,\qquad n_1=- in_2,\qquad l_1=il_2,\qquad l_0=-l_3\;.$$

\underline{Variant $ (23) $}
$$G=\left(\ba{cccc}k_0+k_3&k_1-i k_2&\;\; 2 n_0&0\\
    k_1+i k_2&k_0-k_3&\;\; 2n_1&0\\0&0&\;\; 0&0\\
    2l_1&2l_0&\;\; 2m_1&0\ea\right)$$
$$m_0 =0,\qquad m_3=0,\qquad m_1 =im_2\;,$$
$$n_0= n_3,\qquad n_1= in_2,\qquad l_1=il_2,\qquad l_0=-l_3\;.$$

\underline{Variant $(30)$}
$$G=\left(\ba{cccc}0&2k_1&\;\; n_0+n_3&n_1-i n_2\\
    0&2 k_0&\;\; n_1+i n_2&n_0-n_3\\0&2l_1&\;\; 2 m_0&2m_1\\
    0&0&\;\; 0&0\ea\right),$$
$$l_0 =0,\qquad l_3=0,\qquad l_1 =-il_2\;,$$
$$k_0= -k_3,\qquad k_1= -ik_2,\qquad m_1=-im_2,\qquad m_0=m_3\;.$$

\underline{Variant $ (31) $}
$$G=\left(\ba{cccc}2k_0&0&\;\; n_0+n_3&n_1-i n_2\\
    2k_1&0&\;\; n_1+i n_2&n_0-n_3\\2l_0&0&\;\; 2m_0&2m_1\\
    0&0&\;\; 0&0\ea\right),$$
$$l_1 =0,\qquad l_2=0,\qquad l_0 =l_3\;,$$
$$k_0= k_3,\qquad k_1= ik_2,\qquad m_1=-im_2,\qquad m_0=m_3\;.$$

\underline{Variant $ (32) $}
$$G=\left(\ba{cccc}k_0+k_3&k_1-i k_2&\;\; 0&2 n_1\\
 k_1+i k_2&k_0-k_3&\;\; 0&2n_0\\2l_0&2 l_1&\;\; 0&2m_1\\
 0&0&\;\; 0&0\ea\right),$$
 $$m_0 =0,\qquad m_3=0,\qquad m_1 =-im_2\;,$$
 $$l_0= l_3,\qquad l_1= -il_2,\qquad n_1=-in_2,\qquad n_0 =- n_3\;.$$

\underline{Variant $ (33) $}
 $$G=\left(\ba{cccc}k_0+k_3&k_1-i k_2&\;\; 2n_0&0\\
    k_1+i k_2&k_0-k_3&\;\; 2n_1&0\\2l_0&2l_1&\;\; 2 m_0&0\\
    0&0&\;\; 0&0\ea\right),$$
 $$m_1 =0,\qquad m_2=0,\qquad m_0 =m_3\;,$$
 $$l_0= l_3,\qquad l_1= -il_2,\qquad n_1=in_2,\qquad n_0=n_3\;.$$
Let us remind the general form of the matrix $G$
$$G=\left(\ba{rrrr}k_0+k_3&k_1-i k_2&n_0+n_3&n_1-i n_2\\
    k_1+i k_2&k_0-k_3&n_1+i n_2&n_0-n_3\\l_0+l_3&l_1-i l_2&m_0+m_3&m_1-i m_2\\
    l_1+i l_2&l_0-l_3&m_1+i m_2&m_0-m_3\ea\right).$$
\section{Concluding Remarks}
Let us summarize the results. The matrices of the form
$$
G=\left(\ba{rr}k_0+\; {\bf k}\;
 \vec{\sigma}\;\;&n_0+\; {\bf n}\;\vec{\sigma}\\[3mm]
 l_0+\; {\bf l}\;\vec{\sigma}\;\;&m_0+\; {\bf m}\;\vec{\sigma}\ea\right)\;,$$

\noindent are described by their determinant
\cite{Red'kov-Bogush-Tokarevskaya-2007}
\begin{eqnarray}\det\; G= (kk)\; (mm)+(nn)\;(ll)-2\; (kn)\;(ml)\; -\nonumber\\
 2\;(\;-k_0\;{\bf n}+n_0\; {\bf k}+i\;{\bf k}\times {\bf n}\;)\;
 (\; -m_0\;{\bf l}+l_0\;{\bf m} +i\; {\bf m}\times {\bf l}\;)\; ;\end{eqnarray}
the classification of degenerate 4-dimensional matrices ($\det\, G
=0$),
 of rank 1, 2, 3
can be schematically described by
 $$(k)\longrightarrow  7\;,\;\; (m) \longrightarrow 7\;,\;\;\;\; (n) \longrightarrow  4\;,\;\;\;\; (l) \longrightarrow  4\;,$$
 $$(k,m)\longrightarrow  5\;,\qquad (l,n) \longrightarrow  2\;,\qquad (k,n) \longrightarrow  2\;,$$
 $$(k,l)\longrightarrow  2\;,\qquad (n,m) \longrightarrow  2\;,\qquad (m,l)\longrightarrow  2\;,$$
 $$(k,m,n) \longrightarrow  2\;,\qquad (k,m,l) \longrightarrow  2\;,$$
 $$(n,l,k) \longrightarrow  2\;,\qquad (n,l,m) \longrightarrow  2\;.$$

Thus,   we produced many different subsets of real matrices,
mainly with semi-group structure.
 These subset themselves can be of mathematical interest.
  The question of applying them as Mueller's type  remains open.

\section*{Acknowledgment}

The present work was developed under the auspices of Grant
    F14RA-006, within the cooperation framework between Romanian Academy
    and Belarusian Republican Foundation for Fundamental Research.\par

\end{document}